%% file: main.tex
\documentclass[pdflatex,sn-mathphys-num]{sn-jnl}

\usepackage[colorinlistoftodos,bordercolor=orange,backgroundcolor=orange!20,linecolor=orange,textsize=scriptsize]{todonotes}

\usepackage{geometry}
\usepackage{graphicx}%
\usepackage{multirow}%
\usepackage{amsmath,amssymb,amsfonts}%
\usepackage{amsthm}%
\usepackage{mathrsfs}%
\usepackage[title]{appendix}%
\usepackage{xcolor}%
\usepackage{textcomp}%
\usepackage{manyfoot}%
\usepackage{booktabs}%
\usepackage{algorithm}%
\usepackage{algorithmicx}%
\usepackage{algpseudocode}%
\usepackage{listings}%
\usepackage{subcaption}
\usepackage{adjustbox}
\usepackage{xcolor}
\usepackage{pifont} 
\usepackage{makecell}

\newcommand{\cmark}{\textcolor{green!60!black}{\ding{51}}}
\newcommand{\xmark}{\textcolor{red}{\ding{55}}}

\newtheorem{lemma}{Lemma}

\newcommand{\bx}{\boldsymbol{x}}
\newcommand{\bu}{\boldsymbol{u}}

\newcommand{\by}{\boldsymbol{y}}

\newcommand{\bq}{\boldsymbol{q}}
\newcommand{\bp}{\boldsymbol{p}}
\newcommand{\bz}{\boldsymbol{z}}
\newcommand{\R}{\mathbb{R}}
\newcommand{\Ss}{\mathcal{S}}

\DeclareMathOperator{\latentso}{Latent\_solver}
\DeclareMathOperator{\Diag}{Diag}

\theoremstyle{thmstyleone}%
\newtheorem{theorem}{Theorem}
\newtheorem{proposition}[theorem]{Proposition}%

\theoremstyle{thmstyletwo}%

\theoremstyle{thmstylethree}%

\begin{document}

\title[PI-SONet]{PI-SONet: A Physics-Informed Symplectic Operator Network for Real-Time Optimal Control of Multi-Agent Systems}

\author[1]{\fnm{Alan John} \sur{Varghese}}\email{alan\_john\_varghese@brown.edu}
\equalcont{These authors contributed equally to this work.}

\author[2]{\fnm{Shanqing} \sur{Liu}}\email{shanqing\_liu@brown.edu} 
\equalcont{These authors contributed equally to this work.}

\author[3]{\fnm{Paula} \sur{Chen}}\email{paula.x.chen.civ@us.navy.mil}
\equalcont{These authors contributed equally to this work.}

\author[2]{\fnm{Yaochen} \sur{Zhu}}\email{yaochen\_zhu@brown.edu}

\author*[2]{\fnm{J\'er\^ome} \sur{Darbon}}\email{jerome\_darbon@brown.edu}

\author[2]{\fnm{George} \sur{Em Karniadakis}}\email{george\_karniadakis@brown.edu}

\affil[1]{\orgdiv{School of Engineering}, \orgname{Brown University}, \orgaddress{\city{Providence}, \state{Rhode Island}, \postcode{02906}, \country{U.S.A}}}

\affil[2]{\orgdiv{Division of Applied Mathematics}, \orgname{Brown University}, \orgaddress{\city{Providence}, \state{Rhode Island}, \postcode{02906}, \country{U.S.A}}}

\affil[3]{\orgdiv{Naval Air Warfare Center Weapons Division}, \orgaddress{\city{China Lake}, \state{California}, \postcode{93555}, \country{U.S.A}}}



\abstract{ 
Many real-life applications involve controlling high-dimensional multi-agent systems in real-time. Existing optimal control solvers often suffer from the curse-of-dimensionality and require complete rerunning for each new problem setting.  
We target nonconvex, nonlinear problems in 100s of dimensions by introducing PI-SONet (Physics-Informed Symplectic Operator Network), a structure-preserving operator learning framework for solving parameterized families of optimal control problems and their Pontraygin Maximum Principle (PMP) systems. PI-SONet combines a latent right-space solver with a conditional symplectic operator to produce tractable Hamiltonian trajectories in a computationally efficient auxiliary space and transform them back to physical space.  
This decomposition yields a \textit{single} trained operator that approximates the PMP solution map, 
inherently preserves Hamiltonian structure, and generalizes across unseen problem configurations. 
Unlike existing methods, which are fundamentally single-instance solvers, PI-SONet achieves sub-second inferences on new problem instances, equating to up to 10,000x speedup over representative baselines. 
These results suggest that structure-preserving neural operators provide a practical route toward reusable, real-time surrogates for high-dimensional optimal control.  }



\maketitle




Optimal control provides a natural framework for autonomous decision-making for dynamical systems with applications ranging from robotics and self-driving vehicles to air traffic management and multi-agent coordination~\cite{trelat2005controle,lewis2012optimal}. 
In many practical settings, a critical challenge lies not only in formulating the control problem but in solving it repeatedly across changing environments, initial conditions, and agent capabilities. 
This difficulty is especially pronounced in multi-agent path planning~\cite{foderaro2014distributed, robinson2018efficient,kirchner2020hj,gaspard2025monotone}, where the state dimension grows with the number of agents and collision avoidance, obstacles, and heterogeneous dynamics induce highly nonconvex optimization landscapes.

Classical methods approach these problems from two complementary viewpoints. 
Direct methods~\cite{fahroo2002direct,gong2006pseudospectral,boucher2016galerkin,caillau2023algorithmic} discretize the control problem and solve the resulting nonlinear programming problem but typically must be rerun from scratch for each new problem instance and become expensive at scale. 
Indirect methods either leverage Pontryagin's maximum principle (PMP)~\cite{raymond1998pontryagin,raymond1999pontryagin} to replace the original control problem by a Hamiltonian boundary-value system in the state and co-state variables or dynamic programming (DP) \cite{bellman1957dynamic,bardi2008optimal, sethian2003ordered} to characterize optimality through a Hamilton-Jacobi-Bellman (HJB) equation \cite{Crandall1984TwoAO, falcone2014semi}. 
Despite their mathematical appeal, both paradigms remain limited; DP suffers from the curse of dimensionality, whereas PMP-based solvers generally are highly initialization-sensitive and only recover local optima in nonlinear or nonconvex settings. 

The design of efficient numerical methods for high-dimensional optimal control problems remains a central research challenge. 
Some existing solvers include methods based on the Hopf formula for HJ equations~\cite{darbon2016algorithms,lee2021hopf}, representations using max-plus/tropical algebra~\cite{fleming2000max,Mc2007,akian2008max,dower2015max, akian2023adaptive,akian2026tropical} tensor-decomposition methods~\cite{zbMATH07364328,zbMATH07547920}, trajectory-refinement methods~\cite{alla2019efficient,bokanowski2022optimistic,akian2024multilevel}, and deep learning and neural network-based approaches~\cite{Han2016DeepLA,hure2020deep,kang2021algorithms,nakamura2022neural,darbon2023neural,kunisch2023learning,bokanowski2023neural,sperl2025separable,bottcher2022ai}. 
Notably, recent advances in scientific machine learning and hardware compute power have opened new possibilities for high-dimensional scientific computing~\cite{karniadakis2021physics,cuomo2022scientific,toscano2025pinns,lee2024automatic,lee2025automatic}. In particular, physics-informed neural networks (PINNs) have emerged as an attractive framework for efficiently solving differential equations and constrained dynamical systems in high dimensions~\cite{raissi2019physics,karniadakis2021physics, vinuesa2022enhancing, raabe2023accelerating}. In the context of optimal control, SympOCNet~\cite{meng2022sympocnet} and TSympOCNet~\cite{zhang2025time} show that symplectic neural architectures can solve multi-agent path-planning problems with linear and nonlinear dynamics in 512D. 
These approaches demonstrate  remarkable promise for Hamiltonian-aware learning but remain fundamentally single-instance methods; each new control problem still requires re-solving and retraining for the new corresponding  boundary-value system.

\begin{figure}[!ht]
    \centering
    \includegraphics[width=\linewidth]{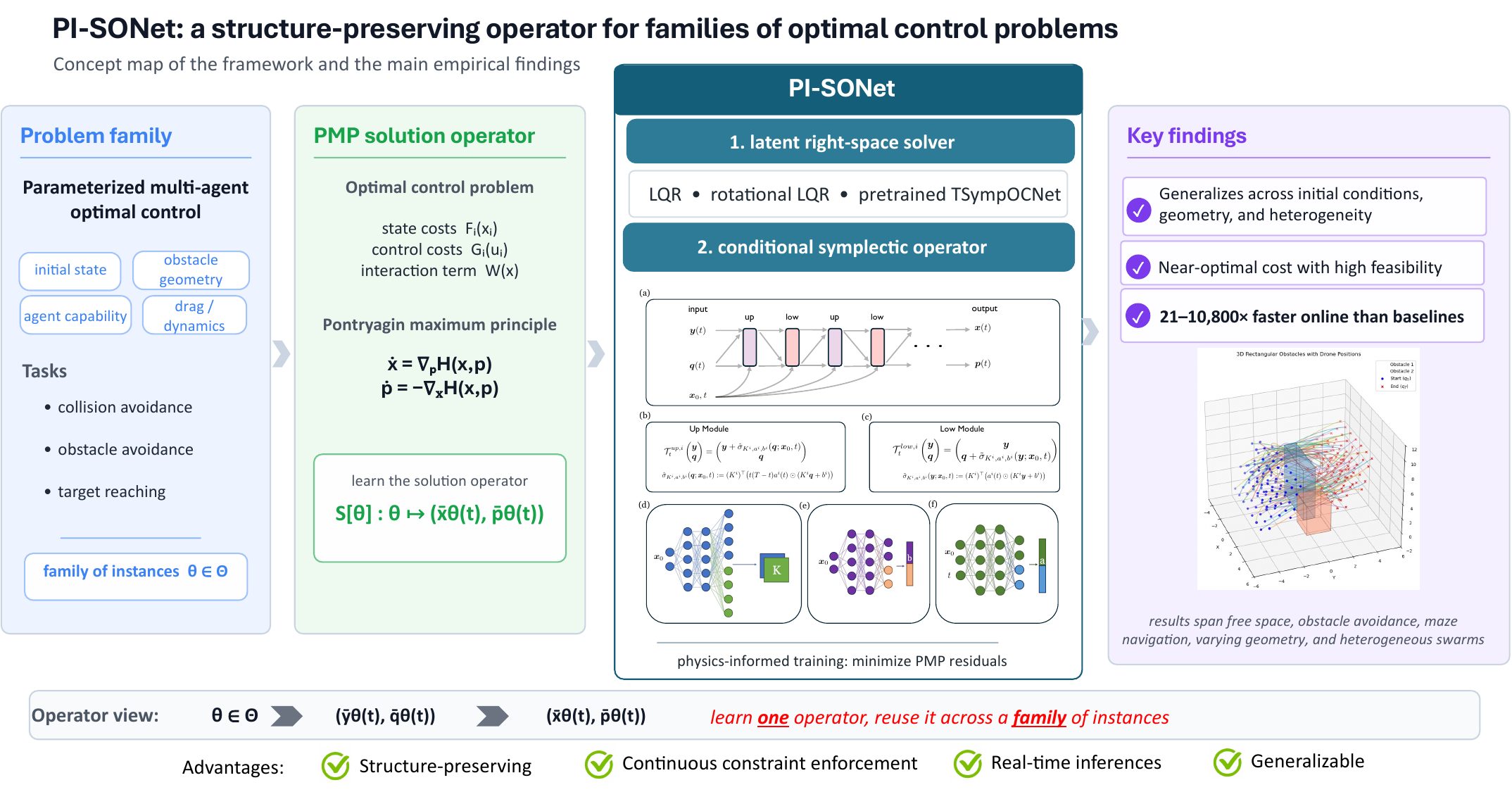}
    \caption{\textbf{PI-SONet: a physics-informed, structure-preserving operator network for families of optimal control problems.}
Parameterization of a multi-agent optimal control problem (e.g., by initial conditions or other problem parameters) induces a corresponding family of Pontryagin maximum principle (PMP) systems and solution operators, rather a than single-instance solution. PI-SONet approximates this family through (1) a latent right-space solver, which produces a computationally tractable Hamiltonian trajectory, and (2) a conditional symplectic operator trained by minimizing PMP residuals that maps the latent trajectory back to physical phase space. On a broad suite of benchmark problems, PI-SONet demonstrates robust generalizability across initial conditions, obstacle geometries, and agent heterogeneity, while achieving near-optimal costs, high feasibility, sub-second online inferences, and $21-10,868x$ speedups over state-of-the-art baseline methods.}
    \label{fig:fig1}
\end{figure}

Instead, we introduce \emph{Physics-Informed Symplectic Operator Network (PI-SONet)}, a \textit{structure-preserving operator-learning framework} for \textit{parameterized families of PMP systems}. Like SympOCNet/TSympOCNet, PI-SONet preserves symplecticity to maintain optimality and feasibility. 
Unlike SympOCNet/TSympOCNet, PI-SONet learns not a single trajectory, but a \emph{family} of solutions over a parameterized class of control problems. 
Neural operators~\cite{lu2021learning,li2020fourier,wang2021learning, brunton2024promising, yin2024scalable,cao2024laplace} are especially appealing in this setting because they 
operate at the level of \textit{function spaces}, learning the solution operator that maps whole families of problem instances to the corresponding families of trajectories, 
rather than approximating individual trajectories in isolation. 
For optimal control, such an approach must address two requirements simultaneously: (1)  preserve the Hamiltonian structure induced by optimality conditions and (2) generalize across meaningful variations of the problem, including changes in geometry, agent properties, and dynamics.  
PI-SONet accomplishes both by combining a latent right-space solver, which produces a tractable Hamiltonian trajectory in a computationally efficient auxiliary space, and a conditional symplectic operator trained by minimizing PMP residuals that maps the latent trajectory back to physical phase space. 
This decomposition yields a robust approximation of the PMP solution operator; a single trained PI-SONet model can be reused across a family of problem instances, while preserving Hamiltonian structure by construction. Figure~\ref{fig:fig1} summarizes this approach.

We view this work as a step toward joining classical numerical optimal control with real-time high-dimensional optimal control. 
Rather than fully replacing existing methods, PI-SONet re-purposes them as latent solvers for simplified surrogate problems and uses a neural operator to map the resulting latent trajectories back to physical phase space. 
In doing so, we extend structured and efficient classical solvers to nonlinear, nonconvex, and/or variable environments (in which classical methods are generally severely limited), while also retaining their numerical and mathematical properties in the latent space.

We evaluate PI-SONet on multi-agent path-planning problems of increasing complexity, including free-space coordination, obstacle avoidance, maze-like environments, variable obstacle geometries, and heterogeneous agent settings in both 2D and 3D physical space. 
Across these scenarios, PI-SONet \textit{generalizes across unseen initial conditions and parameter shifts}, while maintaining \textit{high feasibility }and \textit{near-optimal costs}. 
In benchmark comparisons against multiple shooting, pseudospectral optimal control, and a symplectic PINN baseline, PI-SONet yields substantially faster, \textit{sub-second online inferences} with \textit{speedups spanning $1-4$ orders of magnitude}.
These results suggest that structure-preserving operator learning offers a practical route toward real-time optimal control in complex multi-agent environments. 
More broadly, they indicate that combining Hamiltonian inductive biases with neural operators may help bridge the gap between classical optimal control theory and reusable learned surrogates for families of high-dimensional control problems.

The main contributions of this work are threefold. First, we formulate an operator learning approach for parameterized \textit{families} of PMP systems rather than for isolated trajectories. Second, we introduce a factorization that combines a \textit{tractable} latent Hamiltonian solver with a conditional symplectic decoder, thereby \textit{preserving Hamiltonian structure} while enabling \textit{real-time inferences}. Third, our framework demonstrates \textit{significant advantages over representative baseline methods} on a broad suite of multi-agent path-planning problems involving changes in initial conditions, obstacle geometry, and agent heterogeneity. It is worth noting that while we only demonstrate this framework on path-planning problems, it could also be generalized to other types of control problems given proper PMP formulations and latent solvers.

\section*{Results}\label{sec-result}

We evaluate PI-SONet on a comprehensive suite of multi-agent path-planning problems with the goal of testing three properties: generalizability across initial conditions, adaptivity to changes in geometry and agent capabilities, and computational performance relative to representative baseline methods.
The settings considered include obstacle-free traversal tasks, navigation around circular obstacles or through nonconvex maze-like environments, variable obstacle radii, and heterogeneous agent settings in 2D and 3D. 
In all cases, PI-SONet consistently produces smooth, coordinated trajectories in real-time that remain close to feasible low-cost manifolds, while preserving the Hamiltonian structure of the underlying control problem. Additional details for all numerical examples can be found in the \textit{Supplementary Materials}.

\subsection*{PI-SONet as a structure-preserving operator for families of PMP trajectories}
We begin by viewing PI-SONet not as a solver for a single boundary-value problem but as a learned operator over a \textit{family} of parameterized PMP systems. 
Consider a multi-agent system that consists of $N$ subsystems evolving in  
$\R^{d_x}$ 
over a time horizon $T$. 
For each subsystem $i\in\{1,\dots, N \}$, the state $\bx_i:[0,T]\to \R^{d_x}$ and control $\bu_i:[0,T]\to \R^{d_u}$ satisfy the ODE 
\begin{subequations}\label{dynamics}
\begin{equation}
    \dot{\bx_i} =f_i (\bx_i(s),\bu_i(s)), \ s \in [0,T] 
\end{equation}
with initial condition
\begin{equation}
    \bx_i(0) = \bx_i^0. 
\end{equation}
\end{subequations}
Here, 
$f_i:\R^{d_x} \times \R^{d_u} \to \R^{d_x}$ is a (potentially nonlinear) dynamical map. 

The goal is to steer the system toward a prescribed target, 
encoded by terminal costs $\phi_i^T : \R^{d_x} \to \R$,  while accounting for each subsystem's potential energy $F_i:\R^{d_x}\to\R$, kinetic/control energy $G_i:\R^{d_u}\to\R$, and the interaction energy $W:\R^{N\times d_x}\to\R$ between subsystems and with the environment (e.g., collision and obstacle avoidance in path planning). 
We do so by solving the following problem:
\begin{equation}\label{ocproblem}
    \inf_{\bu(\cdot)}\int_0^T \left( \sum_{i=1}^N \large( F_i(\bx_i(s)) + G_i(\bu_i(s)) \large) + W(\bx(s)) \right) ds + \sum_{i=1}^N \phi^T_i(\bx_i(T))
\end{equation}
subject to~\eqref{dynamics}. 
Here, we denote by $\bx = (\bx_1,\dots,\bx_N)$ and $\bu = (\bu_1,\dots,\bu_N)$ the concatenated states and controls, respectively. 

A necessary optimality condition is then given by PMP as follows. 
Given a costate $\bp = (\bp_1,\dots,\bp_N):[0,T] \to \R^{N \times d_x}$, an optimal state-costate pair $(\bar \bx(\cdot), \bar \bp(\cdot))$ solves the Hamiltonian ODE system
\begin{subequations}\label{hamiltonian_system}
    \begin{equation}
        \left\{ 
        \begin{aligned}
            & \dot{\bar\bx} (s) = \nabla_{\bp} H(\bar \bx(s), \bar \bp(s)), \\ 
            & \dot{\bar \bp} (s) = -\nabla_{\bx} H(\bar \bx(s), \bar \bp(s)) 
        \end{aligned}
        \right.
    \end{equation}
for every time $s \in [0,T]$ and with initial/final boundary condition 
\begin{equation}
    \bar \bx(0) = \bx^0, \ \bar \bp(T) = \nabla_{\bx} \phi^T (\bar \bx (T)), 
\end{equation}
where the Hamiltonian is defined as
\begin{equation}
    H(\bx, \bp ): = \max_{\bu} \left\{ \sum_{i=1}^N \large( \langle \bp_i, f_i(\bx_i, \bu_i) \rangle - (F_i(\bx_i) + G_i (\bu_i) )\large) - W(\bx) \right\}.
\end{equation}

\end{subequations}

For a fixed terminal condition $\phi^T$, the optimal trajectory map associated with~\eqref{ocproblem} 
depends on the initial condition $\bx^0$, the dynamics $f: = (f_1, \dots,f_N)$, the Lagrangian $\ell:=(\ell_1,\dots, \ell_N)$, where  $\ell_i (\bx_i, \bu_i) = F_i(\bx_i) + G_i(\bu_i)$, and the interaction term $W$. 
PI-SONet aims to learn the corresponding solution operator that maps these problem specifications to the PMP solution:
\begin{equation}\label{operator_solution}
    (\bx^0, f, \ell,W) \mapsto \Ss[\bx^0,f,\ell,W](\cdot) = (\bar \bx(\cdot), \bar \bp(\cdot)), \ \text{where $(\bar \bx(\cdot), \bar \bp(\cdot))$ solves~\eqref{hamiltonian_system}} .
\end{equation}
In this work, we focus on \textit{families} of instances for which the problem configuration is parameterized by
$\theta\in\Theta$. 
In our experiments, the parameter $\theta$ collects the initial condition $\bx^0$ together with a low-dimensional set of environment and dynamics parameters (e.g., target configuration, obstacle geometry, agent radius, and drag), while keeping the functional form of~\eqref{ocproblem} fixed.

The workflow for PI-SONet is shown in Figure~\ref{fig:fig1}. As discussed above,  
the main computational challenges arise from nonlinearity in the dynamics, nonconvexity of the
objective/constraints, and high dimensionality due to interactions in multi-agent systems. 
To address these difficulties, PI-SONet proceeds in two stages. 
Note that we emphasize the dependence on the parameterized problem data $\theta$ by writing the PMP solution as 
$(\bar \bx [\theta](\cdot), \bar \bp[\theta](\cdot))$.

\paragraph{Stage 1 (Latent Right-Space Solver).}
We first construct an auxiliary latent problem that is computationally efficient to solve (e.g., using convexification, linearization, or low-dimensional decoupling via sequential subsystem solves). 
A tailored latent solver produces a latent state--costate trajectory
$(\bar \by[\theta](\cdot),\bar \bq[\theta](\cdot))$, i.e., 
\begin{equation}
    \theta \mapsto \latentso [\theta] (\cdot) = ( \bar \by[\theta](\cdot), \bar \bq[\theta] (\cdot)) \ \text{solves a latent Hamiltonian system} . 
\end{equation}
\paragraph{Stage 2 (Structure-Preserving Decoder).}
We construct a structure-preserving symplectic neural operator denoted by $\varphi[\theta,t]$ that depends on the input data $\theta$ and time $t$. 
$\varphi[\theta,t]$ then maps the latent trajectory back to the physical one:
\begin{equation}
    (\bar \by [\theta](t), \bar \bq[\theta](t)) \mapsto \varphi [\theta,t] ( \bar \by [\theta](t), \bar \bq[\theta](t) ) =  ( \bar \bx [\theta](t), \bar \bp[\theta](t) ), \ \forall \ t \in [0,T].
\end{equation}
The neural operator is trained in a physics-informed manner by minimizing residuals of the PMP system~\eqref{hamiltonian_system} to enforce optimality. 
Altogether, the composition of the latent solver and the symplectic operator yields an approximation of the solution
operator in~\eqref{operator_solution}. For a family of inputs $\Theta$,
\begin{equation}
    \varphi[\theta,\cdot] \circ \latentso[\theta] \approx S[\theta], \ \forall \ \theta \in \Theta 
\end{equation}
results in a flexible surrogate that generalizes across problem instances.  
This decomposition enables generalizable solutions of~\eqref{ocproblem} that preserve the structure of the corresponding Hamiltonian system to maintain optimality. Once trained, PI-SONet produces real-time trajectories with favorable scaling in the number of agents.

\subsection*{Multi-agent path planning with collision and obstacle avoidance.} 

\begin{figure}[!ht]
    \centering
    \includegraphics[width=\linewidth]{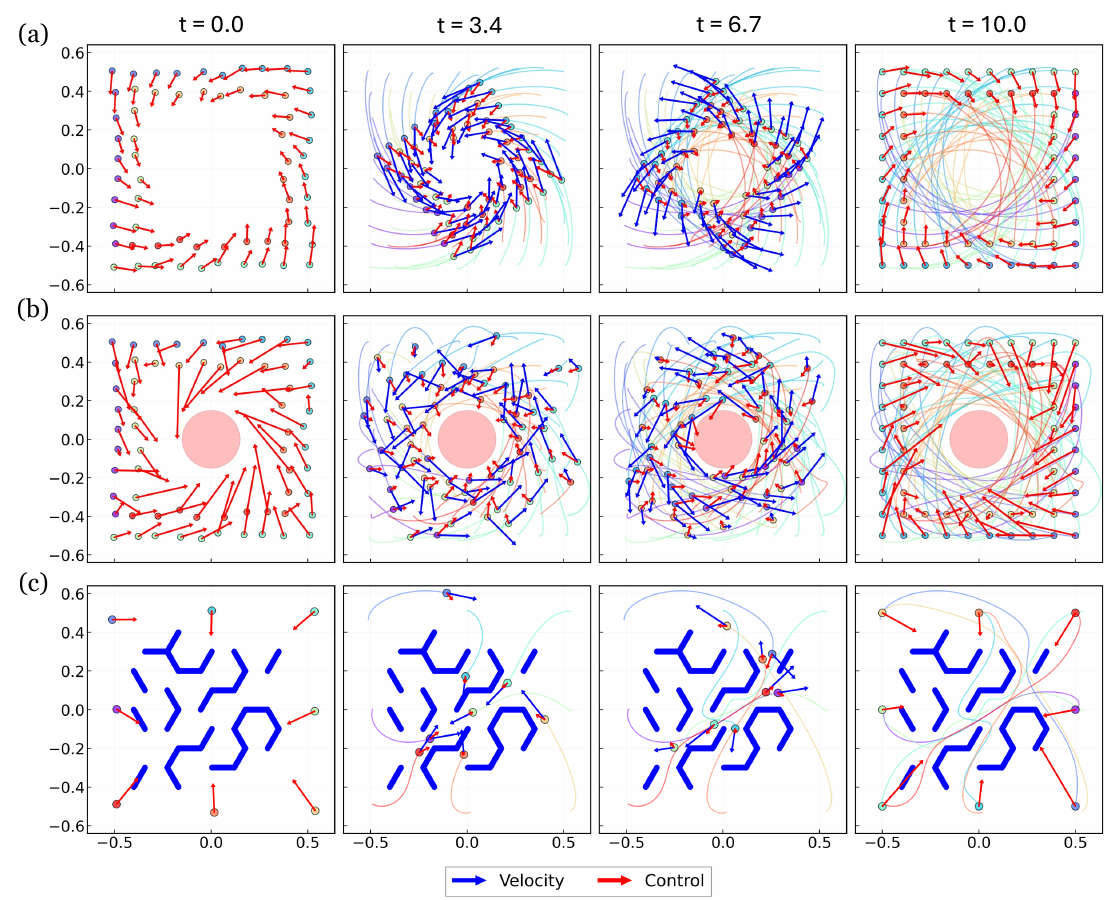}
    \caption{\textbf{Generalization across initial positions in environments of various complexities.}
    Representative PI-SONet predictions at four time instances for three planar multi-agent path-planning tasks: \textbf{(a)} collision avoidance for 64 agents in free space ($256$D overall), \textbf{(b)} navigation of 64 agents around a circular obstacle ($256$D overall), and \textbf{(c)} navigation of 8 agents through a maze-like environment ($64$D overall). Colored curves denote agent trajectories, \textcolor{blue}{blue} arrows denote velocity vectors, and \textcolor{red}{red} arrows denote control inputs. In all three scenarios, PI-SONet produces collision-free coordinated maneuvers in sub-second level timing that adapt to the geometry of the feasible set and are robust to perturbation of the agents' initial positions. These results show that PI-SONet achieves high generalizability, optimality, and feasibility, even in increasingly complex (e.g., constrained, nonconvex/nonsmooth) environments. \textit{See \url{https://github.com/alanjohnvarghese/physics-informed-sonet/tree/main/animations/figure2} for animated visualizations.}}
    \label{fig:multiagent_path_planning}
\end{figure}

We first assess whether a single trained operator can generalize across unseen initial conditions in environments of increasing geometric complexity to ensure that PI-SONet can learn nontrivial solution operators. The goal in each test case is for the agents to swap positions by traversing to the opposite edge/corner of a square domain, while avoiding collisions with each other and obstacles (when present). All agents follow Newtonian dynamics with quadratic-regularized drag ($\dot{\boldsymbol{w}}_i = \boldsymbol{v}_i, \dot{\boldsymbol{v}}_i = \boldsymbol{u}_i - k\boldsymbol{v}_i \|\boldsymbol{v}_i\|_2$ $\forall i = 1, \dots, N$, where $\bx_i =(\boldsymbol{w}_i, \boldsymbol{v}_i)\in\R^4$ denotes the position and velocity of the center of each agent and $k$ is the drag coefficient), 
and the Lagrangian consists of quadratic costs on the velocity and control effort ($\ell_i(\bx_i,\bu_i) = \|\boldsymbol{v}_i\|^2_2 + \|\bu_i\|_2^2$). Each training/test sample then represents the result of a random perturbation of the initial positions $\bx_i$. 
Figure~\ref{fig:multiagent_path_planning} shows representative predictions in three regimes: an obstacle-free space, a domain containing a central circular obstacle, and a maze-like environment. In all three cases, PI-SONet returns coherent position, velocity, and control trajectories and resolves the path-planning task through coordinated collective motion rather than degenerate straight-line crossings.  

The quantitative results show that this behavior remains robust well beyond the visual examples in Figure~\ref{fig:multiagent_path_planning}. In all three settings, PI-SONet achieves 100\% collision-free success on the training set across all tested swarm sizes.
In obstacle-free space, collision-free success on unseen perturbations remains high even at large scale, reaching 96\% for $N=32$, 99\% for $N=56$, and 100\% for $N=64$ agents. 
In the circular-obstacle setting, the corresponding generalization rates remain strong across all tested sizes, including 100\% for $N=4$, 93\% for $N=8$, 99\% for $N=16$, 96\% for $N=32$, 99\% for $N=56$, and 98\% for $N=64$. 
Finally, in the maze setting, which is the most geometrically-restrictive and exhibits the strongest multimodality of the three scenarios, the operator remains extremely robust on unseen perturbations with success rates of 95\% for $N=4$ and 100\% for $N=8$. 

Crucially, these results are achieved with consistently low inference times across all swarm sizes. On an NVIDIA H100 GPU, the forward-pass latency remains at approximately $\sim 0.02$ seconds per run as the number of agents increases from $N=4$ to $N=64$, demonstrating strong scalability in practice. These inference times are compatible with real-time trajectory generation, in contrast to traditional optimal control solvers whose computational cost typically increases drastically with system size.

Qualitatively, we observe that the geometry of the feasible sets strongly shape the learned maneuvers. In free space, the trajectories form braid-like exchanges that maintain pairwise separation while agents swap positions. Introducing a central obstacle induces a distributed rotation around the blocked region, indicating that the operator adapts naturally to convex state constraints. In the maze setting, the predicted trajectories align with the corridor structure and tightly negotiate narrow passageways while preserving inter-agent coordination, demonstrating optimality and feasibility even in nonconvex, nonsmooth environments. 

Note that for any colliding sample, a lightweight test-time refinement step (e.g., 2 steps of LBFGS, which typically takes $< 3$s/sample) that further minimizes the PMP residual can be applied to the failed instance to recover collision-free trajectories. Since the pretrained prediction is typically close to feasible, generally only a small number of extra gradient steps is needed in practice. Unless stated otherwise, the results shown in this paper represent those obtained without this additional fine tuning.



\subsection*{Generalizability to various obstacle geometries}

\begin{figure}[!ht]
    \centering
    \includegraphics[width=\linewidth]{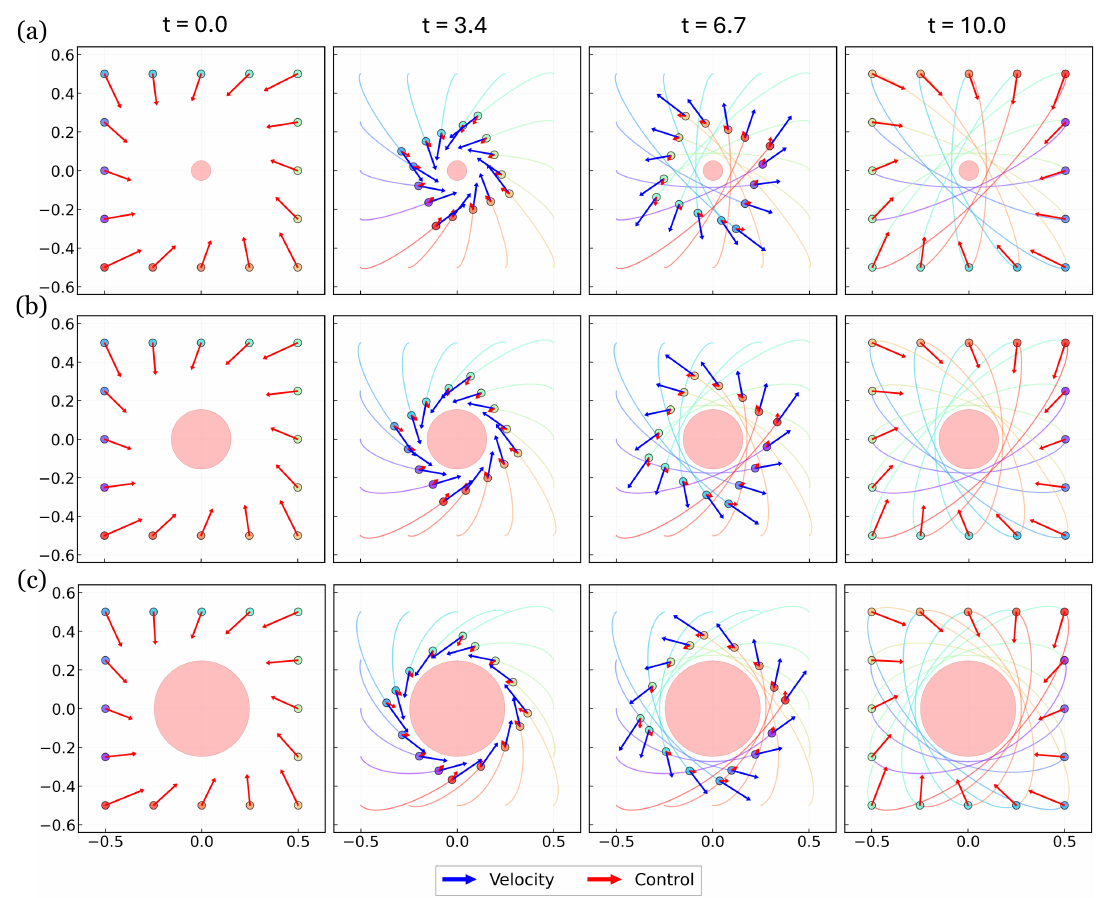}
    \caption{\textbf{Generalization to unseen obstacle geometries.}
    Representative predictions for three different obstacle radii: \textbf{(a)} small, \textbf{(b)} intermediate, and \textbf{(c)} large. Each column corresponds to a different time instant. \textcolor{blue}{Blue} arrows denote velocities, and \textcolor{red}{red} arrows denote controls. As the obstacle radius increases, PI-SONet smoothly expands the trajectories to preserve minimal clearance of the central obstacle, while maintaining the same overall global maneuver. These results demonstrate that PI-SONet is capable of learning \textit{families} of solution operators that continuously depend on abstract problem parameters, such as obstacle geometry. \textit{See \url{https://github.com/alanjohnvarghese/physics-informed-sonet/tree/main/animations/figure3} for animated visualizations.}}
    \label{fig:obstacle_geometry}
\end{figure}

A central advantage of operator learning is that the same model can be conditioned not only on the initial state but also on more abstract problem features, such as environment parameters. Doing so opens possibilities for deploying PI-SONet in non-fixed environments. 
To examine this capability, 
we augment the neural operator input with the obstacle radius and train PI-SONet on a family of circular obstacle problems with varying geometry; i.e., the problem setup remains identical to that of the circular-obstacle case in the previous section but now perturbs the obstacle radius instead of the initial positions. Figure~\ref{fig:obstacle_geometry} shows representative predictions for three different obstacle sizes. 
As the radius of the obstacle grows, the learned trajectories deform smoothly, moving the swarm further from the domain center, automatically increasing obstacle clearance, and overall maintaining organized, collision-free maneuvers. 

The quantitative results confirm that this behavior remains robust and extends beyond the visual examples in Figure~\ref{fig:obstacle_geometry}. PI-SONet achieves 100\% collision-free success on the training set for all tested swarm sizes ($N=4, 8, 16$). Generalization to unseen obstacle geometries remains perfect for $N=4$ and $N=8$ with $50/50$ successful test runs in both cases and remains high for $N=16$ with $48/50$ successful test runs. Together, these results suggest that the learned operator captures a stable and transferable dependence of the PMP solution on obstacle geometry, while preserving feasibility across both training and unseen test instances.

More broadly, the results demonstrate that PI-SONet does indeed approximate a family of solution operators over varying domains, rather than simply memorizing a single obstacle configuration. Notably, the predicted trajectories remain close to the obstacle boundary instead of introducing unnecessarily conservative detours, which indicates that optimality is preserved. This behavior is consistent with the expected trade-off between path length, control effort, and safety constraints inherent to the family of optimal control problems considered in this case and suggests that the learned symplectic operator successfully captures a continuous dependence of the PMP solution on the geometry of the environment.

\subsection*{Handling heterogeneous agents} 
We next consider heterogeneous swarms, in which agent-specific capability parameters are appended to the operator input. By generalizing to agent features, the corresponding operator can be trivially transferred between swarms of perturbed agent capabilities without retraining, potentially allowing for scenarios in which agent capabilities may slowly degrade over time and/or where slightly different agent types are swapped in and out of the swarm. 

\begin{figure}[!ht]
    \centering
    \includegraphics[width=\linewidth]{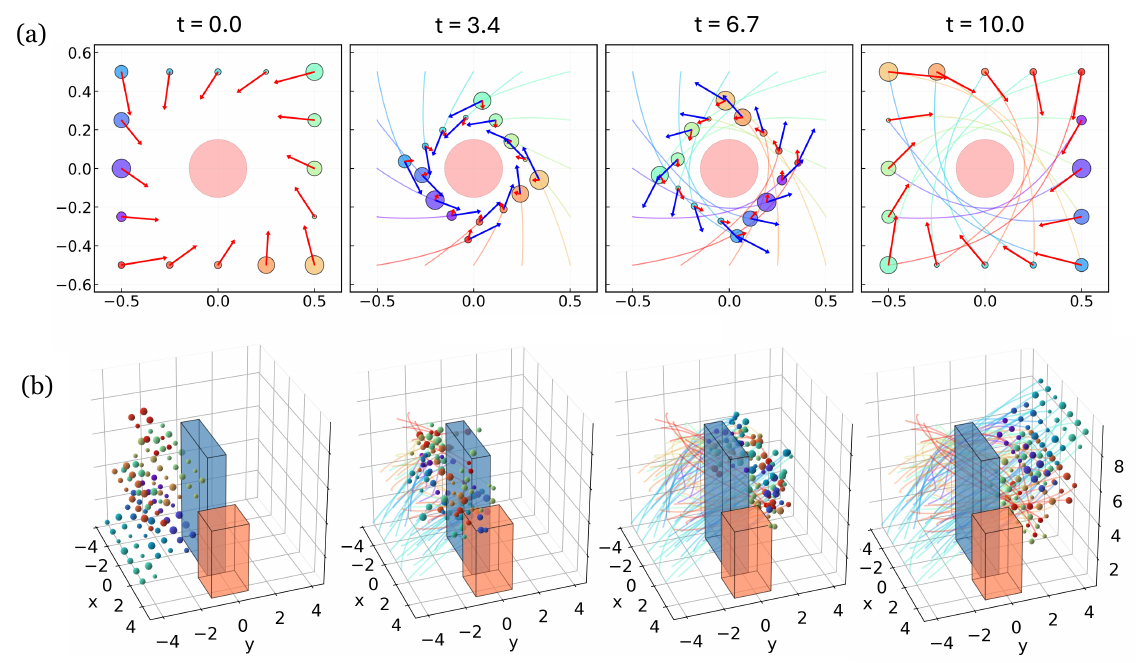}
    \caption{\textbf{PI-SONet handles heterogeneous agents in two and three dimensions.}
    \textbf{(a)} Planar navigation around a circular obstacle with agents of different radii, indicated by marker size. Larger agents induce larger clearance margins while the collective maneuver remains coordinated. \textbf{(b)} Extension to a three-dimensional cluttered environment, where heterogeneous agents navigate around multiple obstacles while preserving spatial separation. These results indicate that PI-SONet can condition on agent-specific capability parameters and transfer across heterogeneous swarms. \textit{See \url{https://github.com/alanjohnvarghese/physics-informed-sonet/tree/main/animations/figure4} for animated visualizations.}}
    \label{fig:heterogeneous_agents}
\end{figure}

Figure~\ref{fig:heterogeneous_agents}\textbf{a} depicts a 2D planar example, in which each agent has a different physical radius (and drag coefficient) and therefore different effective safety margins to model agents of different sizes. Between runs, we perturb the radii (and drag coefficient, proportionally) of all of the agents, but otherwise the problem setup is similar to that in the previous sections.
Again, PI-SONet adapts the coordinated motion accordingly; larger agents execute wider turns and maintain larger clearances from both the obstacle and neighboring agents, whereas smaller agents exploit tighter regions of the domain. 
Importantly, not only does PI-SONet generalize to a perturbation of agent radii, but the agent heterogeneity is also able to be handled by a single shared operator, rather than by training separate models for each agent type.  
Since we also vary the resistance coefficient in this case, we also observe that the learned operator adjusts the control field in a physically-coherent manner; increased drag leads to stronger alignment between control and velocity and hence to a higher total objective cost, reflecting the increased effort required to maintain motion for larger, heavier agents.

We now extend to a much higher dimensional problem by increasing both the physical dimension and swarm size. Figure~\ref{fig:heterogeneous_agents}\textbf{b} shows a cluttered environment in which 100 heterogeneous agents coordinate to traverse a 3D domain with rectangular obstacles. Similarly to before, agents follow Newtonian dynamics but now with much higher dimension ($dx = 6, N = 100$ for overall problem dimension of $600$D) and without drag to simplify the problem slightly. Despite the higher dimensionality, PI-SONet still generates trajectories that preserve spatial separation between agents and adapt obstacle clearance margins to agent size. 
These results suggest that the representation learned by PI-SONet can be extended to higher dimensions, while still maintaining optimality and feasibility. Inferences are again obtained in sub-second level timing (0.0199s for a batch of 50 runs on an NVIDIA H100 GPU), indicating PI-SONet's potential for real-time even in extremely high-dimensional cases.

The quantitative results further support this ability to generalize across heterogeneous settings. In the two-dimensional case with $N=16$ agents, PI-SONet achieves $50/50$ collision-free runs on the training set and retains high performance on unseen test instances with $44/50$ successful runs. In the three-dimensional heterogeneous setting, PI-SONet again achieves $50/50$ collision-free runs during training and generalizes well to unseen test cases with $40/50$ successful runs. In both cases, additional fine tuning at test-time is able to recover 100\% collision-free success rates. These results indicate that 
PI-SONet incorporates both geometric and dynamical heterogeneity within a unified physics-informed operator framework 
and that these properties effectively transfer to very high-dimensional settings.


\subsection*{Benchmarking against state-of-the-art baselines}

We benchmark PI-SONet against three representative baseline numerical solvers: multiple shooting as implemented by CasADi~\cite{SWcasadi,andersson2019casadi}, pseudospectral optimal control as implemented by GPOPS-II~\cite{SWgpops,patterson2014gpops}, and TSympOCNet~\cite{zhang2025time}, which serves as the structure-preserving PINN-analogue to PI-SONet. 
We compare the performance of the methods on the three scenarios from Figure~\ref{fig:multiagent_path_planning}, which exhibit increasingly complex geometries:
\emph{Free}, in which agents cross obstacl-e while avoiding pairwise collisions; \emph{Obstacle}, in which a circular obstacle is added at the center of the domain; and \emph{Maze}, in which agents must navigate through a nonconvex environment. 

The comparison is interpreted in the regime for which PI-SONet is designed. 
The offline training cost is paid once, after which a new problem instance is solved by evaluating the latent solver and conditional symplectic decoder. In other words, the price PI-SONet pays to solve a new instance is the cost of inference, whereas each baseline must be rerun from scratch for each new initial condition. 
Notably, even though TSympOCNet is also symplecticity-preserving, it does not provide generalizability across problem parameters like PI-SONet; TSympOCNet must recommit its entire training time to solve new problem instances, whereas PI-SONet only needs to commit its inference time. 
Under this interpretation, PI-SONet preserves near-optimal costs and high feasibility while dramatically reducing online solution times. 
In Table~\ref{tab:benchmark-sota}, we observe that, across all three settings, PI-SONet achieves speedups of $21-10,868\times$ over the baseline. These gains become more pronounced as the swarm size increases, which demonstrates that whereas most baseline methods suffer from the curse of dimensionality, PI-SONet significantly mitigates it.

The cost and safety statistics follow similar trends. PI-SONet maintains comparable to significantly improved optimality and collision-free rates over the baseline, despite also being generalizable to a continuous range of initial conditions. 
These results indicate that the learned operator achieves large speedups without sacrificing solution quality. We note that while the feasibility of the PI-SONet inferences could be further improved via test-time refinement, the results in Table~\ref{tab:benchmark-sota} reflect PI-SONet's baseline performance without this additional fine tuning.
To further isolate the effect of the symplectic structure, we conduct an ablation study comparing PI-SONet with a standard MLP of comparable capacity. The results show that while feasibility remains similar, the trajectories from MLP fail to adapt to the changing obstacle radius and is suboptimal, highlighting the importance of enforcing symplectic structure (see \textit{Supplementary Materials}).

We also note that direct solvers like multiple shooting and pseudospectral method are highly sensitive to symmetry, initialization, and geometry changes. This sensitivity becomes particularly evident in the \textit{Maze} scenario (see \textit{Supplementary Materials} for qualitative results). 
Moreover, as discretization-based methods, multiple shooting and pseudospectral method only enforce constraints at discrete points and hence, are more susceptible to collisions in complex environments, where coarser grids often miss short-range interactions. 
By contrast, PI-SONet retains the ability to produce feasible warm starts via proper training and its physics-informed architecture 
and, by learning the solution operator, continuously enforces constraints for all time and space, as evidenced by its ability to retain high collision-free rates in the face of highly congested environments and nonconvex obstacles.



\begin{table*}[http]
\centering
\small
\setlength{\tabcolsep}{4pt}
\renewcommand{\arraystretch}{0.95}
\caption{\textbf{Benchmarking PI-SONet against state-of-the-art baselines.}
\textbf{(a)} Average runtime. \textbf{(b)} Average optimal control cost and average maximum collision-constraint violation with the percentage of collision-free runs shown in parentheses. Timing was obtained on a workstation with an Intel Xeon E5-2620 v4 CPU and an NVIDIA GeForce RTX 2080 Ti GPU. The test scenarios are identical to those in Figure~\ref{fig:multiagent_path_planning}. We ran 10 trials for \textit{Free} and \textit{Obstacle} with 4-16 agents, and 5 trials for all other cases. Pseudospectral method has no listed results for the 32-agent cases as it failed to complete after $>12$ hours. For PI-SONet, we report offline training times and online inference times, separately. Speedup is computed as the minimum runtime among the baselines divided by PI-SONet inference time to reflect the time each method requires to re-compute the solution on a new problem instance. PI-SONet demonstrates remarkable speedups, while maintaining high optimality and feasibility, indicating that we do not sacrifice computational speed for solution quality.}
\label{tab:benchmark-sota}

\begin{subtable}[t]{\textwidth}
\centering
\caption{Average runtime.}
\label{tab:benchmark-sota-runtime}
\begin{adjustbox}{width=\textwidth}
\begin{tabular}{ccccccccc}
\toprule
\multirow{2}{*}{\textbf{Scenario}} & \multirow{2}{*}{\textbf{\# Agents}} & \textbf{Multiple} & \multirow{2}{*}{\textbf{Pseudospectral}} & \multirow{2}{*}{\textbf{TSympOCNet}} & \multicolumn{2}{c}{\textbf{PI-SONet}} & \textbf{Minimum} \\
& & \textbf{Shooting} & & & \textbf{Training} & \textbf{Inference} & \textbf{Speedup} \\
\midrule
\multirow{4}{*}{Free}
& 4  & 0.95500s & 11.41528s & 327.69095s & 178.48240s &\textbf{ 0.04492s} & 21.26x \\
& 8 & 9.31030s & 174.51554s & 370.95381s & 254.07590s & \textbf{0.06277s} & 148.33x
 \\
& 16 & 67.22610s & 947.98315s & 475.56235s & 361.48958s & \textbf{0.07472s} & 899.73x \\
& 32 & 2337.97320s & -- & 362.40065s & 493.94911s & \textbf{0.14603s} &	2,481.70x \\
\midrule
\multirow{4}{*}{Obstacle}
& 4  & 2.31660s & 11.39535s & 370.01680s & 230.40773s & \textbf{0.05750s} & 40.29x \\
& 8 & 15.33830s & 302.53436s & 381.33150s & 213.86865s & \textbf{0.06579s }& 233.15x \\
& 16 & 345.11010s & 340.40902s & 397.60065s & 287.69220s & \textbf{0.07493s} & 4,543.26x \\
& 32 & 2341.39220s & -- & 407.28486s & 509.30823s & \textbf{0.13224s} & 3,079.82x \\
\midrule
\multirow{2}{*}{Maze}
& 4  & 27.95880s & 102.40390s & 737.71980s & 419.07924s & \textbf{0.06958s} & 401.85x \\
& 8  & 1518.38500s & 1426.24904s & 900.23110s & 623.88450s & \textbf{0.08283s} & 10,868.39x\\
\bottomrule
\end{tabular}
\end{adjustbox}
\end{subtable}

\vspace{0.75em}

\begin{subtable}[t]{\textwidth}
\centering
\caption{Average cost and safety metrics.}
\label{tab:benchmark-sota-quality}
\begin{adjustbox}{width=\textwidth}
\begin{tabular}{cccccccc}
\toprule
\multirow{2}{*}{\textbf{Scenario}} & \multirow{2}{*}{\textbf{\# Agents}} & \multirow{2}{*}{\textbf{Metric}} & \textbf{Multiple} & \multirow{2}{*}{\textbf{Pseudospectral}} & \multirow{2}{*}{\textbf{TSympOCNet}} & \multicolumn{2}{c}{\textbf{PI-SONet}} \\
& & & \textbf{Shooting} & & & \textbf{Training} & \textbf{Inference} \\
\midrule
\multirow{8}{*}{Free}
& \multirow{2}{*}{4}
& Cost   & 5.1306E-01 & 5.1240E-01 & \textbf{4.9868E-01} & 5.1113E-01 & 5.1168E-01 \\
&        & Safety & 2.00E-03 (90\%) & \textbf{0.00E+00 (100\%)} & \textbf{0.00E+00 (100\%)} & \textbf{0.00E+00 (100\%)} & 1.45E-05 (99\%) \\
\cmidrule(lr){2-8}
& \multirow{2}{*}{8}
& Cost   & 7.7308E-01 & 7.7036E-01 & \textbf{7.5478E-01} & 8.2586E-01 & 8.2750E-01 \\
&        & Safety & 1.40E-02 (0\%) & 1.88E-03 (0\%) &\textbf{ 0.00E+00 (100\%)} & \textbf{0.00E+00 (100\%)} & 1.33E-04 (95\%) \\
\cmidrule(lr){2-8}
& \multirow{2}{*}{16}
& Cost   & 5.2266E+00 & 1.4330E+00         & \textbf{1.4208E+00} & 1.6139E+00 & 1.6112E+00 \\
&        & Safety & 5.49E-03 (10\%) & 5.48E-03 (0\%)   & 1.28E-04 (40\%) & \textbf{0.00E+00 (100\%)} & 6.86E-04 (83\%) \\
\cmidrule(lr){2-8}
& \multirow{2}{*}{32}
& Cost   & \textbf{1.5706E+01} & --         & 2.9847E+00 & 3.3128E+00 & 3.3147E+00 \\
&        & Safety & 1.19E-02 (0\%)  & --         & 3.25E-04 (40\%) & \textbf{0.00E+00 (100\%)} & 1.18E-04 (96\%) \\
\midrule
\multirow{8}{*}{Obstacle}
& \multirow{2}{*}{4}
& Cost   & 5.5662E-01 & 5.5081E-01 & \textbf{5.4929E-01} & 5.7617E-01 & 5.7757E-01 \\
&        & Safety & 2.12E-02 (0\%) & \textbf{0.00E+00 (100\%)} & 3.40E-04 (90\%) & \textbf{0.00E+00 (100\%)} & \textbf{0.00E+00 (100\%)} \\
\cmidrule(lr){2-8}
& \multirow{2}{*}{8}
& Cost   & 8.5410E-01 & 8.5443E-01
 & \textbf{8.3255E-01} & 8.8601E-01 & 8.8959E-01 \\
&        & Safety & 1.32E-02 (0\%) & 2.96E-02 (0\%) & \textbf{0.00E+00 (100\%)} & \textbf{0.00E+00 (100\%)} & 1.99E-04 (93\%) \\
\cmidrule(lr){2-8}
& \multirow{2}{*}{16}
& Cost   & 3.9924E+00 &   1.5544E+00    & \textbf{1.5489E+00} & 1.6244E+00 & 1.6218E+00 \\
&        & Safety & 2.31E-02 (10\%) & 3.36E-05 (0\%)  & \textbf{0.00E+00 (100\%)} & \textbf{0.00E+00 (100\%)} & 2.63E-05 (99\%) \\
\cmidrule(lr){2-8}
& \multirow{2}{*}{32}
& Cost   & 5.0350E+00        & --         & \textbf{3.4541E+00} & 3.8021E+00 & 3.8031E+00 \\
&        & Safety & 1.45E-02 (0\%) & --  & \textbf{0.00E+00 (100\%)} & \textbf{0.00E+00 (100\%)} & 8.06E-05 (96\%) \\
\midrule
\multirow{4}{*}{Maze}
& \multirow{2}{*}{4}
& Cost   & \textbf{2.0744E+00} & 7.4210E+00 & 6.4158E-01 & 2.2645E-01 & 2.2646E-01 \\
&        & Safety & 3.60E-02 (0\%) & 3.44E-02 (0\%) & \textbf{0.00E+00 (100\%)} & \textbf{0.00E+00 (100\%)} & 6.20E-05 (95\%) \\
\cmidrule(lr){2-8}
& \multirow{2}{*}{8}
& Cost   & 1.6805E+01 & 1.8489E+01 & \textbf{1.4654E+00} & 2.4217E-01 & 2.4379E-01
 \\
&        & Safety & 3.49E-02 (0\%) & 3.48E-02 (0\%) & 2.71E-02 (20\%)
 & \textbf{0.00E+00 (100\%)} & \textbf{0.00E+00 (100\%)} \\
\bottomrule
\end{tabular}
\end{adjustbox}
\end{subtable}
\end{table*}

\section*{Discussion}\label{sec12}

We introduced PI-SONet, a physics-informed, structure-preserving operator network for solving families of optimal control problems. Our extensive numerical experiments show that PI-SONet provides a robust operator-level approximation of families of PMP trajectories across a broad range of multi-agent control settings. 
Namely, a single trained model is able to generalize across unseen initial conditions, varying obstacle geometries, and heterogeneous agent capabilities, while retaining high feasibility, near-optimal costs, and real-time inferences, outperforming representative state-of-the-art baselines on nearly every metric. 
Taken together, our numerical experiments indicate that structure-preserving operator learning moves beyond single-instance trajectory prediction by providing reusable real-time surrogates for parameterized optimal control problems.
In particular, the learned operator remains effective even when the geometry or dynamics of the problem varies within a low-dimensional family, as reflected by its smooth adaptation to obstacle radius, agent size, and drag in both planar and three-dimensional settings. These results suggest that the conditional symplectic representation successfully captures 
the underlying dependence of PMP trajectories on problem parameters.

While our 3D example demonstrates impressive scaling to very high dimensions (600D overall) with millisecond inference times, deployment on real-life systems may require a decentralized approach to maintain scalability on resource-constrained hardware platforms. 
As such, one natural future direction is to reformulate PI-SONet using a divide-and-conquer approach, in which each agent (or each small subsystem of agents) is associated with an individual operator, while an additional global operator enforces coordination and consistency across the full swarm. In addition to improving scalability, such a hierarchical decomposition would potentially also reduce computational complexity, improve system modularity, and better exploit the partial decoupling structure that is often present in large multi-agent systems.

Another direction is to combine the present PMP-based operator framework with HJ/DP-based methods. In particular, one may solve a low-dimensional HJB equation for each agent (or each local subsystem) to obtain a feedback control or value-based prior and then use this information to construct a stronger latent solver for the global operator. In contrast to PMP solvers, which only guarantee local optimality, HJ/DP approaches generally yield globally optimal solutions, even in the presence of nonconvex Hamiltonians. As a result, incorporating HJ-based latent solvers may provide more principled right-space representations for adversarial or competitive settings, in which nonconvex Hamiltonians are involved and PMP-based solvers may be inadequate. 

Overall, the PI-SONet framework combines classical optimal control with the computational and approximation power of neural operators. 
In doing so, we move beyond the current state-of-the-art for optimal control, which is otherwise generally limited in dimensionality, problem complexity, and generalizability to problem settings. 
Instead, PI-SONet offers promising pathways toward real-time optimal control in complex multi-agent environments. Namely, by extending the architecture and/or latent right-space solvers in PI-SONet, our approach provides methodological routes for further improvements in scalability, robustness, and physical interpretability.


\section*{Methods}\label{sec-method}

We seek to learn, in operator form, the family of optimal trajectories associated with a parameterized class of multi-agent optimal control problems. 
The central idea is to factor the solution operator into two components: a fast \emph{right-space} solver that produces a latent Hamiltonian trajectory and a \emph{conditional symplectic neural operator} that decodes this latent trajectory back to physical phase space. 
This factorization separates the geometric and nonconvex features of the physical problem from the computational burden of repeatedly solving a high-dimensional boundary-value problem.

\subsection*{Parameterized optimal control problem and Pontryagin system}

In this section, we define the solution operator for the parameterized families of optimal control problems that PI-SONet aims to learn. Let $\vartheta\in\Theta$ denote the problem instance. 
In the applications we consider, $\vartheta$ collects the initial condition together with a low-dimensional set of environment and dynamics parameters, such as target configuration, obstacle geometry, agent radii, and drag coefficients. 
For $N$ agents with states $\bx_i:[0,T]\to\R^{d_x}$ and controls $\bu_i:[0,T]\to\R^{d_u}$, $i = 1, \dots,N$, we consider an optimal control problem of the following form: 
\begin{subequations}\label{eq:method_ocp}
\begin{align}
\min_{\bu(\cdot)}\quad 
& \int_0^T \Bigg[
\sum_{i=1}^N\Big(F_i^\vartheta(\bx_i(t))+G_i^\vartheta(\bu_i(t))\Big)
+ U_{\epsilon,\ell}\!\big(h_\vartheta(\bx(t))\big)
\Bigg]\,dt
+ \phi_\vartheta(\bx(T)), \label{eq:method_ocp_cost}\\
\text{s.t.}\quad 
& \dot \bx_i(t)=f_i^\vartheta(\bx_i(t))+B_i^\vartheta \bu_i(t), \qquad i=1,\dots,N,\label{eq:method_ocp_dyn}\\
& \bx(0)=\bx_0^\vartheta, 
\label{eq:method_ocp_ic}
\end{align}
\end{subequations}
where $\bx=(\bx_1,\dots,\bx_N)$, $u=(\bu_1,\dots,\bu_N)$, $B_i^\vartheta\in\R^{d_x\times d_u}$, $F_i^\vartheta$ is an individual state cost, $G_i^\vartheta$ is a convex individual control penalty, $U_{\epsilon,\ell}$ is a smooth barrier/penalty function used to relax hard state constraints, and $h_\vartheta$ encodes obstacle and collision avoidance.  
The terminal cost $\phi_\vartheta$ may represent either a soft target cost or, in the fixed-endpoint case, an exact penalty enforcing $\bx(T)=\bx_T^\vartheta$.

We assume throughout that, for each $\vartheta\in\Theta$, the maps $f_i^\vartheta$, $F_i^\vartheta$, $h_\vartheta$, and $\phi_\vartheta$ are continuously differentiable and that each $G_i^\vartheta$ is proper, differentiable, and strictly convex. We define the 
associated Hamiltonian by
\begin{equation}\label{eq:method_H}
H_\vartheta(\bx,\bp)
=
\sum_{i=1}^N
\left(
\langle \bp_i,f_i^\vartheta(\bx_i)\rangle
-
F_i^\vartheta(\bx_i)
+
\big(G_i^\vartheta\big)^*\!\big((B_i^\vartheta)^\top \bp_i\big)
\right)
-
U_{\epsilon,\ell}\!\big(h_\vartheta(\bx)\big),
\end{equation}
where $(G_i^\vartheta)^*$ denotes the Legendre--Fenchel transform of $G_i^\vartheta$. 
PMP associated with the above optimal control problem introduces a costate $\bp$ that lies in the dual space of the state to yield the following state-costate system:
\begin{subequations}\label{eq:method_pmp}
\begin{align}
\dot \bx(t) &= \nabla_{\bp} H_\vartheta(\bx(t),\bp(t)), \label{eq:method_pmp_x}\\
\dot \bp(t) &= -\nabla_{\bx} H_\vartheta(\bx(t),\bp(t)), \label{eq:method_pmp_p}\\
\bx(0) &= \bx_0^\vartheta,\qquad
\bp(T)=\nabla \phi_\vartheta(\bx(T)). \label{eq:method_pmp_bc}
\end{align}
\end{subequations}

Let us denote by
\[
\mathcal S:\Theta\to C([0,T];\R^{2Nd_x}),
\qquad
\mathcal S[\vartheta](t)=(\bar \bx_\vartheta(t),\bar \bp_\vartheta(t)),
\]
the solution operator that maps the problem instance $\vartheta$ to the corresponding state-costate trajectory. The following lemma allows us to reconstruct the physical control once the co-state has been predicted. 
\begin{lemma}\label{lem:control_recovery}
Assume that $G_i^\vartheta$ is differentiable and strictly convex for every $i$ and $\vartheta$. Then, the maximizer in \eqref{eq:method_H} is unique and is given pointwise by
\begin{equation}\label{eq:method_control_recovery}
\bar \bu_i(t)=\nabla \big(G_i^\vartheta\big)^*\!\big((B_i^\vartheta)^\top \bar \bp_i(t)\big).
\end{equation}
\end{lemma}

\subsection*{Right-space factorization of the solution operator}
Directly learning the solution operator $\mathcal S$ in physical space is difficult when \eqref{eq:method_pmp} is high-dimensional, nonlinear, and/or nonconvex. 
We instead introduce a latent Hamiltonian system that is cheaper to solve and use it as our right-space representation. To this end, 
 for each instance $\vartheta$, let
  \[
    \mathcal T_{\mathrm{lat}}[\vartheta](t)=(\bar \by_\vartheta(t),\bar \bq_\vartheta(t))
  \]
be the trajectory generated by a latent Hamiltonian solver. The only structural requirement is that $(\bar \by_\vartheta,\bar \bq_\vartheta)$ solves a Hamiltonian ODE
\begin{equation}\label{eq:method_latent_general}
\begin{aligned}
\dot \by(t)&=\nabla_{\bq} \widetilde H_\vartheta(\by(t),\bq(t)) \\ 
\dot \bq(t)&=-\nabla_{\by} \widetilde H_\vartheta(\by(t),\bq(t)), 
\end{aligned}
\end{equation}
for some smooth latent Hamiltonian $\widetilde H_\vartheta$. 

In the numerical examples in the \textbf{Results} section, we primarily employ linear-quadratic regulator (LQR)-type latent solvers.
 A representative choice for an LQR-type solver is
\begin{equation}\label{eq:method_latent_H}
\widetilde H_\vartheta(\by,\bq)
=
\sum_{i=1}^N\left(
\langle \bq_i, A_i^\vartheta \by_i\rangle
-\frac12 \by_i^\top Q_i^\vartheta \by_i
+\frac12 \bq_i^\top C_i^\vartheta \bq_i
\right),
\end{equation}
where $A_i^\vartheta\in\R^{d_x\times d_x}$, $Q_i^\vartheta\in\R^{d_x\times d_x}$ and $C_i^\vartheta\in\R^{d_x\times d_x}$ are chosen to provide a tractable prior. Such choices include standard LQR-type latents and topology-aware variants obtained by adding skew-symmetric rotational couplings. In the standard linear-quadratic case, each agent evolves according to
\begin{equation}\label{eq:method_latent_linear_system}
\frac{d}{dt}
\begin{pmatrix}
\by_i\\ \bq_i
\end{pmatrix}
=
\underbrace{
\begin{pmatrix}
A_i^\vartheta & C_i^\vartheta\\
Q_i^\vartheta & - (A_i^\vartheta)^\top
\end{pmatrix}
}_{=: \mathcal H_i^\vartheta}
\begin{pmatrix}
\by_i\\ \bq_i
\end{pmatrix},
\end{equation}
which can be solved efficiently using matrix exponentials or other symplectic ODE solvers. 
The full solution operator is then approximated as the composition
\begin{equation}\label{eq:method_operator_factorization}
\widehat{\mathcal S}_\omega[\vartheta](t)
=
\Phi_\omega^\vartheta\!\left(t,\bar \by_\vartheta(t),\bar \bq_\vartheta(t)\right),
\end{equation}
where $\Phi_\omega^\vartheta$ is a time-dependent symplectic decoder parameterized by trainable weights $\omega$ for each fixed $\vartheta$. 

Note that this right-space factorization is not limited to the LQR setting. 
More generally, the latent problem may be obtained by simplifying the original optimal control problem using any means that preserves computational tractability. 
For example, removing collision-avoidance or other interaction terms yields a separable system of $N$ independent low-dimensional optimal control problems, one per agent. 
These reduced problems can then be solved efficiently using classical methods, including direct methods, shooting methods, or grid-based HJ solvers in moderate dimensions. 
In this sense, the latent solver is best viewed as a flexible numerical prior adapted to the structure of the simplified problem.

\subsection*{Conditional symplectic neural operator} 

\begin{figure}
\centering
    \includegraphics[width=\linewidth]{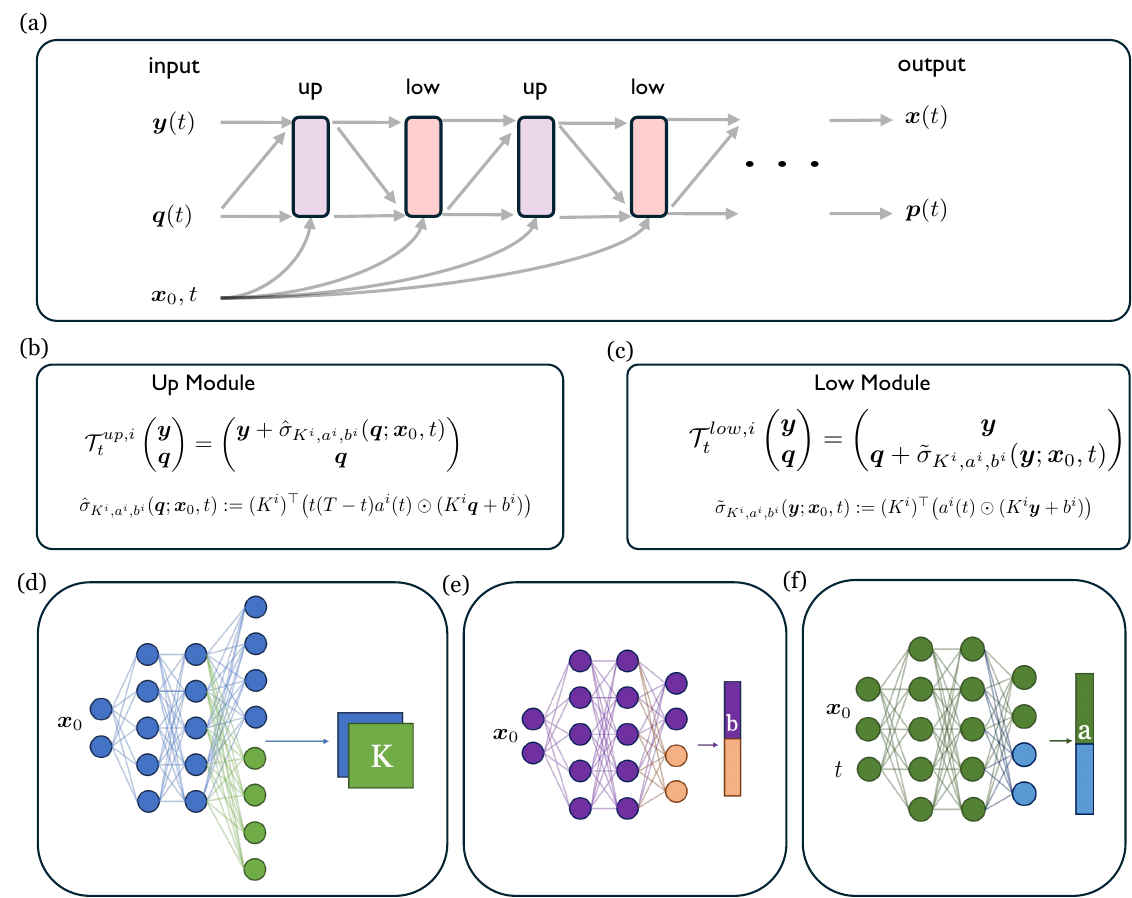}
    \caption{\textbf{Overall architecture of the conditional symplectic decoder Symplectic Operator Network (SONet).}  A latent Hamiltonian trajectory $(\by(t),\bq(t))$ is mapped to the physical trajectory $(\bx(t),\bp(t))$ through alternating lower and upper symplectic blocks, whose parameters are generated from the problem inputs. This architecture yields a time-dependent family of structure-preserving maps that generalizes across problem instances.}
   \label{fig:sonet_architecture}
\end{figure}

\subsubsection*{Symplectic architecture}

Hamiltonian systems, like those arising from PMP, have symplectic structure. In this section, we describe the symplecticity-preserving architecture that underlies our conditional symplectic decoder. Let $\bz=(\by,\bq)\in\R^{2n}$ with $n=Nd_x$, and let
\[
J=
\begin{pmatrix}
0&I_n\\
-I_n&0
\end{pmatrix}.
\]
A differentiable map $\Phi:\R^{2n}\to\R^{2n}$ is symplectic if
\[
(\nabla \Phi(\bz))^\top J\,\nabla \Phi(\bz)=J
\qquad
\text{for all } \bz\in\R^{2n}.
\]

To preserve this structure, we define $\Phi_\omega^\vartheta$ as a composition of conditional triangular symplectic shears. 
For layer $\ell=1,\dots,L$, let
\begin{equation}\label{eq:method_param_nets}
K_\ell(\vartheta)=\mathcal N_{K,\ell}(\vartheta),\qquad
b_\ell(\vartheta)=\mathcal N_{b,\ell}(\vartheta),\qquad
a_\ell(\vartheta,t)=\mathcal N_{a,\ell}(\vartheta,t),
\end{equation}
where $\mathcal N_{K,\ell}$, $\mathcal N_{b,\ell}$, and $\mathcal N_{a,\ell}$ are small neural networks. We then define the scalar potential as
\begin{equation}\label{eq:method_potential}
\Psi_\ell^\vartheta(\bz,t)
=
\frac12 (K_\ell(\vartheta)\bz+b_\ell(\vartheta))^\top
\Diag(a_\ell(\vartheta,t))
(K_\ell(\vartheta)\bz+b_\ell(\vartheta))
\end{equation}
and its gradient as
\begin{equation}\label{eq:method_sigma}
\sigma_\ell^\vartheta(\bz,t)
=
\nabla_{\bz} \Psi_\ell^\vartheta(\bz,t)
=
K_\ell(\vartheta)^\top
\Big(
a_\ell(\vartheta,t)\odot (K_\ell(\vartheta)\bz+b_\ell(\vartheta))
\Big).
\end{equation}
For fixed-endpoint problems, we use the boundary-preserving factor
\[
\beta(t)=t(T-t),
\]
and define the lower and upper shears by
\begin{align}
\mathcal T_{\ell,t}^{\mathrm{low},\vartheta}(\by,\bq)
&=
\big(\by,\; \bq+\sigma_\ell^\vartheta(\by,t)\big), \label{eq:method_low}\\
\mathcal T_{\ell,t}^{\mathrm{up},\vartheta}(\by,\bq)
&=
\big(\by+\beta(t)\sigma_\ell^\vartheta(\bq,t),\; \bq\big). \label{eq:method_up}
\end{align}
The conditional symplectic operator is the alternating composition
\begin{equation}\label{eq:method_decoder}
\Phi_\omega^\vartheta(t,\cdot)
=
\mathcal T_{L,t}^{\mathrm{up},\vartheta}
\circ
\mathcal T_{L,t}^{\mathrm{low},\vartheta}
\circ \cdots \circ
\mathcal T_{1,t}^{\mathrm{up},\vartheta}
\circ
\mathcal T_{1,t}^{\mathrm{low},\vartheta}.
\end{equation}

\begin{proposition}\label{prop:shear_symplectic}
For every instance $\vartheta$, time $t\in[0,T]$, and layer $\ell$, the maps
$\mathcal T_{\ell,t}^{\mathrm{low},\vartheta}$ and $\mathcal T_{\ell,t}^{\mathrm{up},\vartheta}$
defined in \eqref{eq:method_low}--\eqref{eq:method_up} are symplectic.
\end{proposition}

\begin{proposition}\label{prop:decoder_symplectic}
For every $\vartheta\in\Theta$ and $t\in[0,T]$, the decoder $\Phi_\omega^\vartheta(t,\cdot)$ in \eqref{eq:method_decoder} is symplectic. Moreover, if $\beta(0)=\beta(T)=0$, then its state component preserves endpoints in the sense that
\[
\big(\Phi_\omega^\vartheta(0,\by,\bq)\big)^{(\by)}=\by,
\qquad
\big(\Phi_\omega^\vartheta(T,\by,\bq)\big)^{(\by)}=\by.
\]
\end{proposition}

Proposition~\ref{prop:shear_symplectic} follows from the fact that
$\nabla_{\bz}\sigma_\ell^\vartheta(\bz,t)=K_\ell(\vartheta)^\top \Diag(a_\ell(\vartheta,t)) K_\ell(\vartheta)$
is symmetric, so the Jacobians of the shear maps are triangular symplectic matrices. Proposition~\ref{prop:decoder_symplectic} then follows by closure of the symplectic group under composition.

\subsubsection*{Multi-agent parameterization} 

For large swarms, a dense matrix $K_\ell(\vartheta)\in\R^{n\times n}$, $n = Nd_x$ quickly becomes expensive since its output dimension grows quadratically with the number of agents. We therefore impose a block-diagonal parameterization
\begin{equation}\label{eq:method_blockdiag}
K_\ell(\vartheta)=
\Diag\!\big(
K_{\ell,1}(\vartheta),\dots,K_{\ell,N}(\vartheta)
\big),
\qquad
K_{\ell,i}(\vartheta)\in\R^{d_x\times d_x}
\end{equation}
and generate the blocks using a shared trunk with agent-wise output heads. 

\begin{proposition}[Scaling of the block-diagonal decoder]\label{prop:block_scaling}
Let $n=Nd_x$. A dense parameterization of $K_\ell(\vartheta)\in\R^{n\times n}$ requires $\mathcal O(n^2)=\mathcal O(N^2d_x^2)$ outputs per layer, whereas the block-diagonal parameterization \eqref{eq:method_blockdiag} requires only $\mathcal O(Nd_x^2)$ outputs. Thus, for fixed per-agent state dimension, the output dimension of the decoder grows linearly rather than quadratically with the number of agents.
\end{proposition}

This restriction is particularly useful in the operator learning regime, where the same architecture must generalize across many initial conditions and environment parameters. Moreover, this choice of structure balances expressivity of the operator with computational efficiency and scalability. In our numerical experiments, we observe that this relatively sparse structure provides sufficiently high fidelity and model expressivity, leading to no noticeable degradation in the accuracy of the learned solution operator.

\subsection*{Physics-informed training}

We enforce the physics of the PMP system via the training loss for PI-SONet. Given a problem instance $\vartheta$, we compute the latent trajectory
$(\bar \by_\vartheta,\bar \bq_\vartheta)=\mathcal T_{\mathrm{lat}}[\vartheta]$ and define the decoded trajectory by 
\begin{equation}\label{eq:method_decoded}
(\hat \bx_\omega^\vartheta(t),\hat \bp_\omega^\vartheta(t))
=
\Phi_\omega^\vartheta\!\left(t,\bar \by_\vartheta(t),\bar \bq_\vartheta(t)\right).
\end{equation}
Its time derivative is then obtained by the chain rule as follows:
\begin{equation}\label{eq:method_chain_rule}
\frac{d}{dt}
\begin{pmatrix}
\hat \bx_\omega^\vartheta(t)\\
\hat \bp_\omega^\vartheta(t)
\end{pmatrix}
=
D_{\bz}\Phi_\omega^\vartheta\!\left(t,\bar \by_\vartheta(t),\bar \bq_\vartheta(t)\right)
\begin{pmatrix}
\dot{\bar \by}_\vartheta(t)\\
\dot{\bar \bq}_\vartheta(t)
\end{pmatrix}
+
\partial_t \Phi_\omega^\vartheta\!\left(t,\bar \by_\vartheta(t),\bar \bq_\vartheta(t)\right),
\end{equation}
which is evaluated using automatic differentiation.

We train $\omega$ by minimizing the residual of the physical PMP system. For a batch $\mathcal B\subset\Theta$ of sampled instances and collocation times $\{t_j\}_{j=1}^{N_t}$, we define
\begin{align}
r_{\bx}^\vartheta(t_j;\omega)
&=
\dot{\hat \bx}_\omega^\vartheta(t_j)
-
\nabla_{\bp} H_\vartheta\big(\hat \bx_\omega^\vartheta(t_j),\hat \bp_\omega^\vartheta(t_j)\big),
\label{eq:method_res_x}\\
r_{\bp}^\vartheta(t_j;\omega)
&=
\dot{\hat \bp}_\omega^\vartheta(t_j)
+
\nabla_{\bx} H_\vartheta\big(\hat \bx_\omega^\vartheta(t_j),\hat \bp_\omega^\vartheta(t_j)\big).
\label{eq:method_res_p}
\end{align}
The training objective is then given by
\begin{equation}\label{eq:method_loss}
\mathcal L(\omega)
=
\frac{1}{|\mathcal B|N_t}
\sum_{\vartheta\in\mathcal B}\sum_{j=1}^{N_t}
\left(
\|r_{\bx}^\vartheta(t_j;\omega)\|_2^2
+
\|r_{\bp}^\vartheta(t_j;\omega)\|_2^2
\right)
+
\lambda_0 \mathcal L_{\mathrm{IC}}
+
\lambda_T \mathcal L_{\mathrm{TC}},
\end{equation}
where
\begin{align}
\mathcal L_{\mathrm{IC}}
&=
\frac{1}{|\mathcal B|}
\sum_{\vartheta\in\mathcal B}
\|\hat \bx_\omega^\vartheta(0)-\bx_0^\vartheta\|_2^2,\label{eq:method_ic_loss}\\
\mathcal L_{\mathrm{TC}}
&=
\frac{1}{|\mathcal B|}
\sum_{\vartheta\in\mathcal B}
\|\hat \bp_\omega^\vartheta(T)-\nabla \phi_\vartheta(\hat \bx_\omega^\vartheta(T))\|_2^2.
\label{eq:method_tc_loss}
\end{align}
In fixed-endpoint problems, the summand in \eqref{eq:method_tc_loss} is replaced by the terminal state mismatch loss
$\|\hat \bx_\omega^\vartheta(T)-\bx_T^\vartheta\|_2^2$.

\subsection*{Continuation strategy for state constraints}

The barrier term $U_{\epsilon,\ell}$ is introduced to regularize collision and obstacle constraints during training. 
To progressively enforce feasibility, we use a continuation strategy over the barrier parameters:
\begin{equation}\label{eq:method_continuation}
\epsilon_{m+1}=\rho_\epsilon \epsilon_m,\qquad
\ell_{m+1}=\rho_\ell \ell_m,
\qquad 0<\rho_\epsilon,\rho_\ell<1,
\end{equation}
starting from a smooth, weakly penalized problem and warm-starting each stage from the optimizer state of the previous one. 
At stage $m$, the Hamiltonian $H_\vartheta$ in \eqref{eq:method_H} is replaced by $H_{\vartheta,\epsilon_m,\ell_m}$, and the objective \eqref{eq:method_loss} is minimized until convergence. This annealed training procedure reduces optimization stiffness in the early phase of training, while steering the learned trajectories toward the feasible set as the penalty sharpens.

\subsection*{Inference}

After training, inference for a new instance $\vartheta$ requires only two steps: (i) solve the latent Hamiltonian system to obtain $(\bar \by_\vartheta,\bar \bq_\vartheta)$ and (ii) evaluate the decoder $\Phi_\omega^\vartheta$ along that latent trajectory. The resulting computational pathway
\[
\vartheta
\;\mapsto\;
(\bar \by_\vartheta,\bar \bq_\vartheta)
\;\mapsto\;
(\hat \bx_\omega^\vartheta,\hat \bp_\omega^\vartheta)
\]
provides an amortized approximation of the physical PMP solution operator. In particular, the expensive training stage is performed only once, whereas the online cost for a new problem instance reduces to a latent solve and a single forward pass through the symplectic operator, which can typically be done in sub-second level timing.

\backmatter


\section*{Acknowledgments}

This research was developed with funding from the Defense Advanced Research Projects Agency (DARPA). The views, opinions and/or findings expressed are those of the author and should not be interpreted as representing the official views or policies of DARPA or the U.S. Government. Specifically, this research is supported by the DARPA The Right Space (TRS) grant HR00112590073 and the Laboratory University Collaboration Initiative (LUCI) program through an award made by the Office of the Under Secretary of War for Research and Engineering (OUSW(R\&E)), Science and Technology (S\&T)/Foundations.  Distribution Statement A. Approved for Public Release; Distribution is unlimited. PR 26-0070. 

\section*{Competing Interests}

\bibliography{ref,references}

\newpage
\input{Supplementary.tex}

\end{document}

%% file: Supplementary.tex






 \newcommand{\bX}{\boldsymbol{X}}
 \newcommand{\bv}{\boldsymbol{v}}









    

\setcounter{figure}{0}
\setcounter{table}{0}
\setcounter{equation}{0}
\renewcommand{\thefigure}{S\arabic{figure}}
\renewcommand{\thetable}{S\arabic{table}}
\renewcommand{\theequation}{S\arabic{equation}}
 
\begin{center}
{\Large\bfseries Supplementary Materials for ``PI-SONet: A Physics-Informed Symplectic Operator Network for Real-Time Optimal Control of Multi-Agent Systems''\footnote{Distribution Statement A. Approved for Public Release; Distribution is unlimited. PR 26-0070.}}
\end{center}
\vspace{1em}

\section{Implementation Details}\label{secA1}
\subsection{Block diagonal structure of K}

\begin{figure}[!ht]
    \centering
    \includegraphics[width=0.9\linewidth]{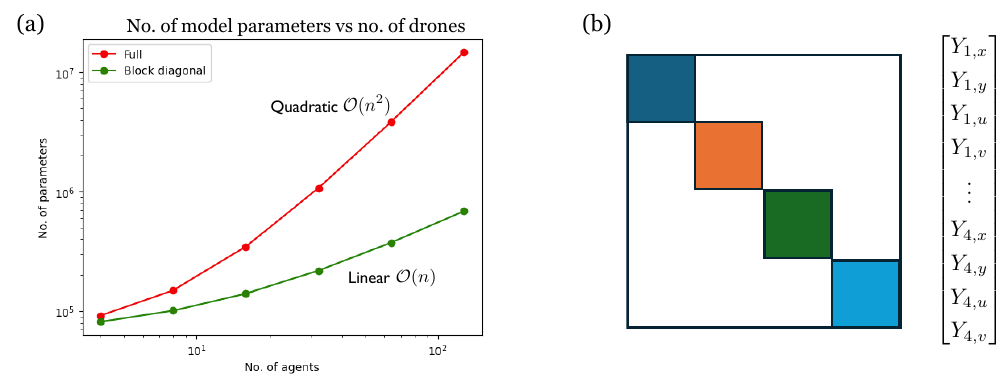}
    \caption{\textbf{Block diagonal structure of $K$.} (a) Scaling of the number of SONet parameters with the number of drones for full vs. block-diagonal parameterizations of $K$. (b) Schematic illustration of the block-diagonal structure we impose for $K$, where each diagonal block corresponds to an individual drone and cross-drone blocks are set to zero. By imposing a relatively sparse structure for $K$, we balance expressivity of SONet with computational scalability. In our numerical experiments, we do not observe any noticeable 
    degradation in accuracy as a result of this sparsity.}
    \label{fig:block_diagonal}
\end{figure}

The matrix $K$ has dimension $Nd_x \times Nd_x$, where $N$ denotes the number of drones and $d_x$ is the dimension of the state space (with $d_x=4$ in our planar experiments and $d_x = 6$ in our 3D experiments). As the number of drones increases, learning a full $K$ matrix quickly becomes impractical since the output layer of the network that predicts $K$ scales quadratically with $n$, resulting in prohibitively larger numbers of parameters. 

To address this issue, we impose a block-diagonal structure on the matrix $K$. We assume that only the blocks corresponding to individual drones are nonzero, while blocks associated with cross-drone terms are set to zero. This parameterization of $K$ matrix reduces the number of model parameters and improves scalability.

In addition, the the network prediction $K$ shares weights and biases across all pre-final layers. The final layer is implemented using multiple heads (as shown in Figure~5(d)), where each head outputs the block of $K$ corresponding to a single drone. This design further improves parameter efficiency while maintaining sufficient representational capacity. Empirically, we observe that enforcing a block-diagonal structure on $K$ does not degrade the performance of the operator network in our experiments.

\begin{table}[h]
\centering
\begin{tabular}{c|cc}
\hline
\textbf{No. of drones} & \textbf{Full} & \textbf{Block diagonal} \\
\hline
4  & 91,848   & 81,480  \\
8  & 149,448  & 101,064 \\
16 & 347,592  & 140,232 \\
32 & 1,075,656 & 218,568 \\
64 & 3,858,888 & 375,240 \\
\hline
\end{tabular}
\caption{\textbf{Number of model parameters for full vs. block-diagonal parameterizations of $K$ as a function of the number of drones.} Full, dense structure quickly leads to infeasibly large numbers of network parameters, whereas the block diagonal structure significantly mitigates this effect.}
\label{tab:parameter_scaling}
\end{table}

Figure
~\ref{fig:block_diagonal} illustrates how the number of parameters in SONet scales with the number of drones for both the full and block-diagonal parameterizations of K; corresponding values are also reported in Table~\ref{tab:parameter_scaling}. When a full $K$ matrix is used, the number of model parameters scales quadratically with the number of drones. In contrast, the block-diagonal structure leads to linear scaling with respect to $n$, resulting in a substantially more parameter-efficient model.

\subsection{Latent solver}
\subsubsection{LQR solver with rotation}
The standard LQR formulation provides an efficient latent solver for linear dynamics but implicitly favors straight-line trajectories. In multi-agent settings with many drones operating in a confined domain, standard LQR induces a degenerate topology; trajectories tend to intersect near the center of the domain, creating a high-density bottleneck. As a result, the downstream operator must significantly deform these intersecting trajectories to satisfy collision constraints, leading to optimization difficulties and poor scalability.

In the space-time domain $[0,T] \times \mathbb{R}^2$, the joint trajectories of the agents form a braid. Feasible, collision-free solutions correspond to non-intersecting strands, i.e., trajectories lying in a non-trivial homotopy class of the braid space \cite{mavrogiannis2020decentralized}. For symmetric initial and terminal configurations, these feasible solutions typically exhibit a rotational structure, whereas the standard LQR solution lies in a degenerate class with strand intersections.

To address this mismatch while retaining the efficiency of the standard LQR solver, we introduce a rotational component into the latent dynamics. Specifically, we augment the system matrix with a skew-symmetric term:
\begin{equation}
    A_{\text{rot}} = \begin{pmatrix} 0 & I \\ 0 & \Omega \end{pmatrix}, \quad 
    \Omega = \begin{pmatrix} 0 & -C_B \\ C_B & 0 \end{pmatrix}.
\end{equation}

This modification induces curvature in the latent trajectories, encouraging coordinated rotational motion, thereby avoiding the central bottleneck. This skew-symmetric structure ensures that this component does not introduce artificial dissipation, preserving the conservative nature of the dynamics.

From a geometric perspective, this rotational bias steers the latent trajectories toward a non-intersecting braid class, providing a more suitable initialization for the downstream operator. Empirically, this behavior leads to significantly improved robustness and scalability compared to standard LQR, particularly in high-density regimes (see Section~\ref{subsec:latent_solver}).

\subsubsection{Eikonal reference solver}
To find a PMP trajectory in a highly non-convex environment one has to initialize the trajectory properly so that the loss converges to a good basin of attraction.
The reference trajectory solves a stochastic particle system via Euler-Maruyama on the SDE
$$d\bX_t = f(\bX_t, t)\,dt + \sigma\, dW_t,$$
where $\bX_t \in \mathbb{R}^{2N}$ stacks the positions of $N$ agents.
Drift is given by $f = G \cdot \bv(\bX)$, where $G \approx I$ is a perturbed-identity Riemannian metric and $\bv = (\bv_1, \dots, \bv_N)$ ($\bv_i\in \mathbb{R}^{2}, \forall i$) is a sum of vector fields acting on each agent $i$ given by
$$\bv_i(X) =\bv_{\text{rep},i}(\bX) + \bv_{\text{wall},i}(\bX) + \bv_{\text{nav},i}(\bX),$$
where inter-agent repulsion (gradient of a log-barrier) is given by
$$\bv_{\text{rep},i}(\bX) = -\sum_{j \neq i} \frac{\bx_i - \bx_j}{c_1\left(\|\bx_i - \bx_j\| - 2r_{\text{agent}}\right)^2},$$
wall repulsion (gradient of an inverse-square barrier from the closest segment of the wall $\pi_k(\bx_i)$) is given by
$$\bv_{\text{wall},i}(\bX) = \sum_{k} \frac{\bx_i - \pi_k(\bx_i)}{c_2\left(\|\bx_i - \pi_k(\bx_i)\| - r_{\text{agent}} - r_{\text{obs}}\right)^3},$$
and the navigation vector field is given by
$$\bv_{\text{nav},i}(\bX) = -\frac{\nabla u(\bx_i)}{\|\nabla u(\bx_i)\|},$$
which is evaluated by bilinear interpolation of the precomputed gradient field of the Eikonal navigation $u: \mathbb{R}^2 \to \mathbb{R}$. We take $u$ to be the solution to the Eikonal equation
$$\|\nabla u(\bx)\| = \frac{1}{s(\bx)}, \quad u(\bx^*) = 0,$$
where $s(\bx) \approx 0$ inside obstacles and $s(\bx) = 1$ in free space. 
Note that $\bv_{\text{nav},i}$ is not the gradient of any closed-form potential. Hence, we compute the gradient as follows.
Five trials are run with different random seeds; the best run is selected by minimizing the detour ratio:
$$\mathcal{C} = \max_i \frac{\int_0^T \|\dot{\bx}_i(t)\|\,dt}{\|\bx_i(0) - \bx_i^*\|}.$$
Next, we rescale the resulting trajectory in time to match the final time of the original optimal control problem. The rescaled reference trajectory is then used to pretrain the SONet by training the loss function $\mathcal{L}(\theta) = \|\bx_\theta - \bx_{ref}\|^2$, where $\bx_\theta$ represents the transformed trajectory. We note that in our experiments we used this only for the \textit{Maze} case.
 
\subsection{Perturbed initial conditions}
\subsubsection{Crossing a square (`free'): scalability to large swarms}

To evaluate the scalability of PI-SONet to high-dimensional state spaces, we consider a multi-agent room-crossing task in the absence of obstacles. The agents must move from perturbed initial positions to antipodal target locations while avoiding inter-agent collisions. The complexity of this problem arises purely from pairwise interactions, which scale quadratically with the number of drones. We consider swarm sizes ranging from $N=4$ to $N=64$ agents.

The SONet architecture consists of 3 up layers and 3 down layers across all cases, with a network width of 8 for both $K$ and $b$. Latent trajectories are generated using an LQR solver with rotation.

To ensure feasibility as the swarm density increases, we adjust the physical parameters across different regimes. For $N \in \{4, 8, 16, 32\}$, the drone radius is fixed at $0.020$ with a perturbation radius of $0.05$. For denser swarms ($N=56$ and $N=64$), the drone radius is reduced to $0.016$, and the perturbation radius is adjusted to $0.04$ and $0.025$, respectively.

Training strategies differ between low- and high-dimensional regimes. For smaller swarms ($N \in \{4, 8, 16\}$), we use constant penalty parameters $\epsilon = 10^{-4}$ and $\lambda = 10^{-4}$ and train using 150 Adam steps followed by 100 L-BFGS steps. For larger swarms ($N \in \{32, 56, 64\}$), we employ an annealing strategy, initializing $\epsilon$ and $\lambda$ at $0.1$ and decaying them by a factor of 0.6 every 20 epochs. In these cases, the model is trained using 100 Adam epochs followed by 100 L-BFGS steps. Additionally, weight decay ($L^2$ regularization) is applied across all experiments with values ranging from $10^{-6}$ to $5 \times 10^{-5}$.

\begin{figure}[!ht]
    \centering
    \includegraphics[width=\linewidth]{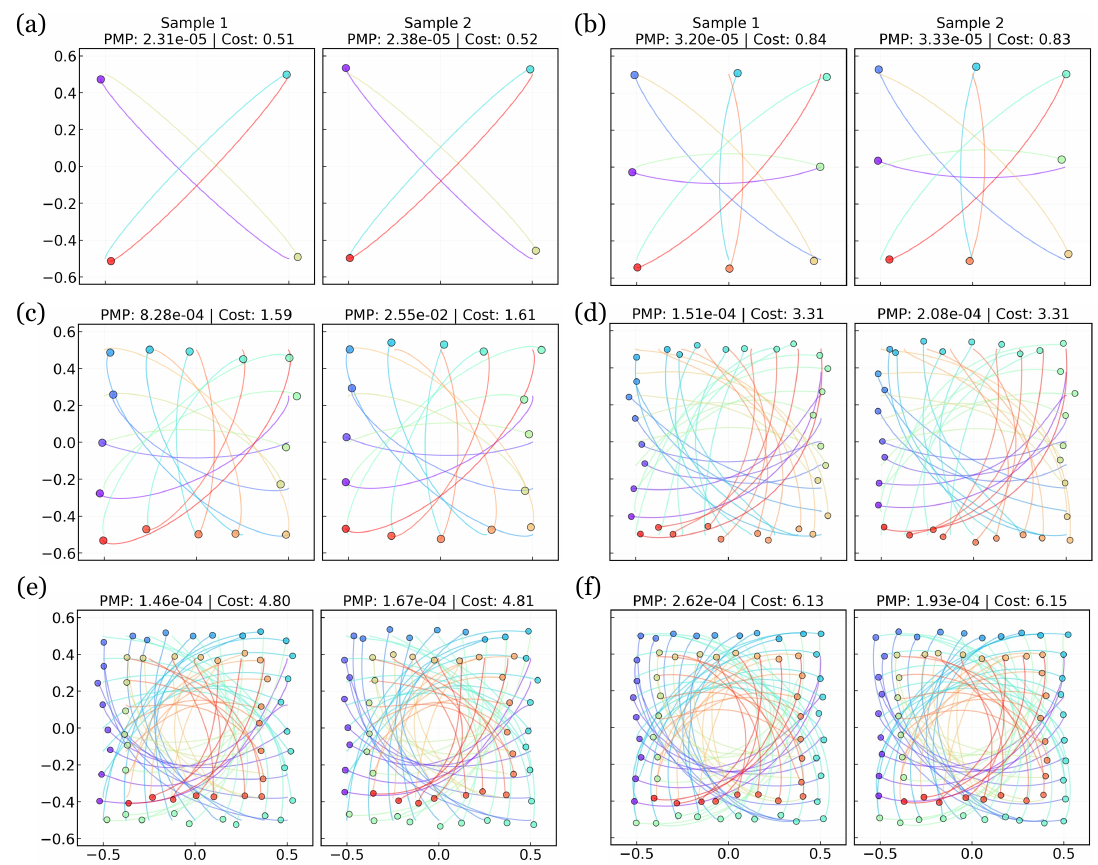}
    \caption{\textbf{Scalability in an obstacle-free room.} Trajectories for swarms of varying sizes ($N \in \{4, \dots, 64\}$) navigating an obstacle-free room. For each swarm size, we show two samples of perturbed initial positions from the test dataset. PI-SONet exhibits high optimality, high feasibility, and real-time inferences across all swarm sizes, demonstrating efficient scalability for high-dimensional trajectory planning problems.}
    \label{fig:swarm_trajectories}
\end{figure}

Figure~\ref{fig:swarm_trajectories} visualizes representative trajectories for two test samples for each swarm size. PI-SONet demonstrates strong scalability, achieving 100\% collision-free success rates on the training set across all swarm sizes. Generalization to unseen test perturbations remains robust: 99\% for 4 agents, 95\% for 8 agents, 83\% for 16 agents, 97\% for 32 agents, 99\% for 56 agents, and 100\% for 64 agents.

In addition to accuracy, PI-SONet exhibits consistently low inference latency across all swarm sizes. On an NVIDIA H100 GPU, the average batched forward-pass time remains on the order of $10^{-2}$ seconds with measured runtimes of $0.018$s ($N=4$), $0.015$s ($N=8$), $0.019$s ($N=16$), $0.021$s ($N=32$), $0.019$s ($N=56$), and $0.019$s ($N=64$). All inference times (in this and all other test scenarios) are computed on a single batch consisting of all test samples. These results indicate that the method maintains efficient inferences even as the number of agents increases, thereby enabling real-time trajectory generation in practice.

\begin{table}[h]
\centering
\resizebox{\textwidth}{!}{%
\begin{tabular}{c|cccccc}
\hline
\textbf{No. of drones} & Drone radius & Perturbation radius & Annealing & Weight decay ($\lambda)$  & $C_B$ & $C_Q$\\
\hline
4  & 0.020  & 0.05 & No & Yes ($1e-6$) & $\pi/20$ & 1.0   \\
8  & 0.020  & 0.05  & No& Yes ($1e-6$ & $\pi/20$ & 0.1\\
16 & 0.020  & 0.05 & No & Yes ($5e-5$)& $\pi/20$ & 0.05 \\
32 & 0.016  & 0.05 & Yes& Yes ($1e-6$)& $\pi/10$ & 0.001 \\
56 & 0.016  & 0.04 & Yes& Yes ($1e-6$)& $\pi/10$ & 0.001\\
64 & 0.016  & 0.025 & Yes& Yes ($1e-6$)& $\pi/10$ & 0.001\\
\hline
\end{tabular}
}
\caption{\textbf{Hyperparameters and physical settings for obstacle-free scalability experiments.} The table details the configuration for swarms of varying sizes $N$ navigating an obstacle-free environment. As the agent density increases ($N \ge 32$), the drone radii and initial perturbation radii are reduced to maintain geometric feasibility. Additionally, training stability is enforced through annealing, weight decay ($\lambda$), and adjusted coefficients for the rotation prior ($C_B$) and quadratic velocity ($C_Q$) cost terms in the latent LQR solver.}
\label{tab:free_table}
\end{table}

\begin{table}[!ht]
\resizebox{\textwidth}{!}{%
\centering
\begin{tabular}{c c c c c c c}
\toprule
\# Drones & Pass (Train) & Pass (Test) & Avg PMP (Train) & Avg PMP (Test) & Avg $L$ (Train) & Avg $L$ (Test) \\
\midrule
4  & 100/100 & 99/100 & $9.02 \times 10^{-5}$ & $1.71 \times 10^{-4}$ & 0.511 & 0.512 \\
8  & 100/100 & 95/100 & $7.93 \times 10^{-5}$ & $1.21 \times 10^{-3}$ & 0.826 & 0.828 \\
16 & 100/100 & 83/100 & $1.14 \times 10^{-3}$ & $7.70 \times 10^{-3}$ & 1.614 & 1.612 \\
32  & 100/100 & 97/100 & $1.33 \times 10^{-4}$ & $7.39 \times 10^{-4}$ & 3.313 & 3.315 \\
56  & 100/100 & 99/100 & $1.38 \times 10^{-4}$ & $2.50 \times 10^{-4}$ & 4.807 & 4.809 \\
64 & 100/100 & 100/100 & $1.91 \times 10^{-4}$ & $2.12 \times 10^{-4}$ & 6.123 & 6.122 \\
\bottomrule
\end{tabular}}
\caption{\textbf{Performance of PI-SONet for the `free' case with perturbed initial position.} PI-SONet exihibits high optimality and feasibility across all swarm sizes tested.}
\label{tab:free_performance_table}
\end{table}

These results suggest that the combination of the symplectic prior and the annealing strategy effectively structures the optimization landscape, such that PI-SONet avoids poor local minima and maintains high performance even in high-dimensional, dense-interaction settings.

\subsubsection{Crossing a square with a circular obstacle (`obstacle'): Scalability to Large Swarms}

We extend the scalability analysis to include a static environmental constraint. In this setting, a swarm of size $N \in \{4, \dots, 64\}$ performs an antipodal position swap while avoiding a fixed circular obstacle of radius $r_{\text{obs}} = 0.15$ located at the center of the domain. The objective remains to reach target states from perturbed initial positions without inter-agent or obstacle collisions.

The presence of the obstacle introduces a nontrivial topological constraint by blocking direct straight-line paths between antipodal pairs. As a result, agents must learn coordinated trajectories that circumnavigate the obstacle. We employ LQR with rotation and pretrained TSympOCNet as the TRS solver to guide trajectory generation.

The SONet architecture consists of 3 up layers and 3 down layers across all cases. The network width for $K$ and $b$ is set to 4 for $N \in \{4,8,16\}$, 16 for $N=32$, and 8 for $N \in \{56,64\}$. 

Training strategies vary with swarm size. For $N \in \{4,8,16\}$, we use constant penalty parameters $\epsilon = \lambda = 10^{-4}$ and train using 150 Adam steps followed by 200 L-BFGS steps. For $N=32$, we increase the penalty parameters to $\epsilon = \lambda = 10^{-3}$ and train using 150 Adam steps followed by 50 L-BFGS steps. In the aforementioned low-dimensional cases ($N \leq 32$), no annealing is applied. 

For larger swarms ($N \in \{56,64\}$), we employ annealing to stabilize optimization in the high-dimensional regime. The penalty parameters are initialized at $\epsilon = \lambda = 0.1$ and decayed by a factor of 0.6 every 20 steps. The model is trained using 100 Adam steps followed by 50 L-BFGS steps for $N=56$ and 100 Adam steps followed by 30 L-BFGS steps for $N=64$.

\begin{figure}[!ht]
    \centering
    \includegraphics[width=\linewidth]{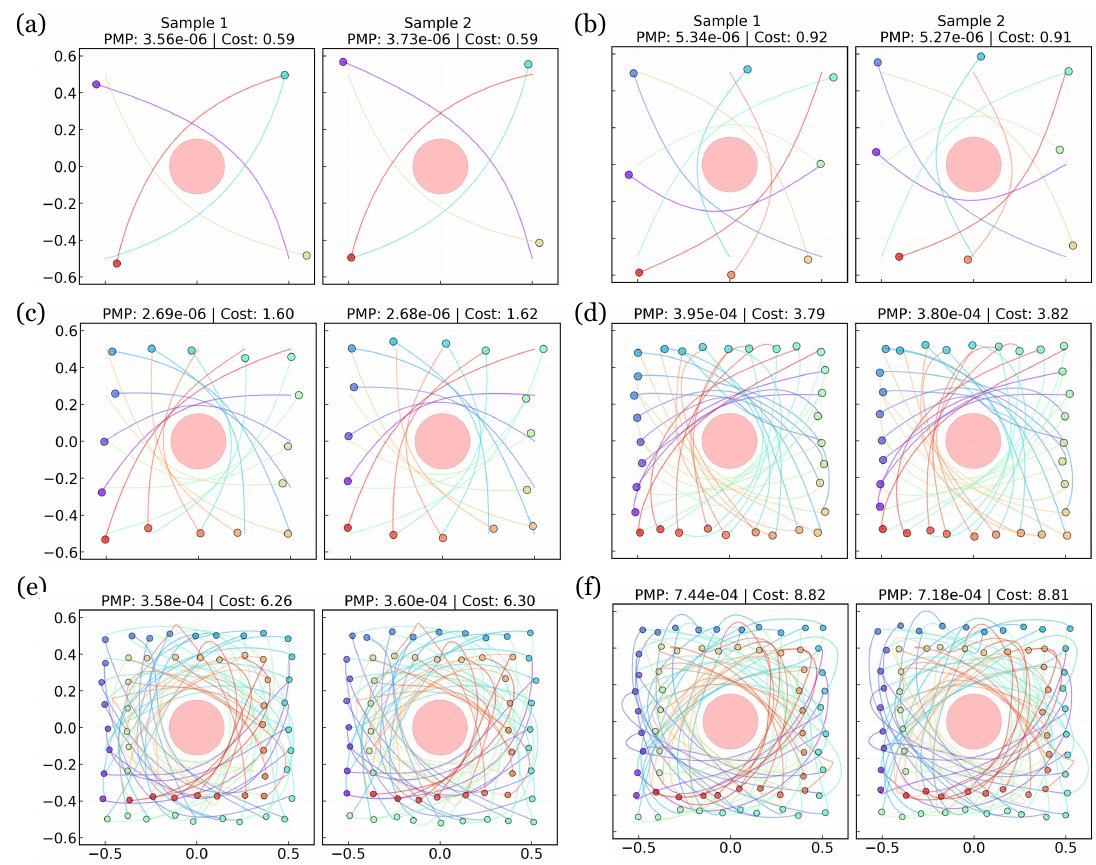}
    \caption{\textbf{Scalability with a central obstacle.} Trajectories for swarms of sizes ($N \in \{4, \dots, 64\}$) navigating around a fixed circular obstacle ($r=0.15$). The obstacle prevents agents from crowding the center of the domain, inducing a coordinated circular flow that simplifies inter-agent avoidance. As such, PI-SONet's performance in this scenario is in some cases better than that in the obstacle-free scenario. For each case, we show two test samples with perturbed initial positions.}\label{fig:sonet_obstacle}
\end{figure}

The inclusion of the obstacle yields consistently strong performance. PI-SONet achieves a 100\% collision-free success rate on the training set across all swarm sizes. Generalization to unseen test perturbations remains robust, with success rates of 100\% (4 agents), 94\% (8 agents), 99\% (16 agents), 97\% (32 agents), 99\% (56 agents), and 98\% (64 agents).

In addition to accuracy, PI-SONet maintains consistently low inference latencies across all swarm sizes in the presence of obstacles. On an NVIDIA H100 GPU, the average forward-pass time per batch remains on the order of $10^{-2}$ seconds with measured runtimes of $0.017$s ($N=4$), $0.017$s ($N=8$), $0.019$s ($N=16$), $0.018$s ($N=32$), $0.019$s ($N=56$), and $0.020$s ($N=64$). These results demonstrate that the method retains efficient inferences even when faced with additional geometric constraints, which shows PI-SONet's potential for real-time trajectory generation in obstacle-rich environments.

Interestingly, performance in the obstructed setting is comparable to, and in some cases slightly better than, the obstacle-free case. This behavior can be attributed to the geometric role of the obstacle; in the absence of obstacles, antipodal trajectories tend to converge toward the origin, creating a high-density bottleneck. The central obstacle eliminates this degeneracy by forcing agents to circumnavigate the center, inducing a coordinated circular flow that naturally distributes agents in space and reduces collision complexity despite the reduced free volume.

Figure~\ref{fig:sonet_obstacle} illustrates representative trajectories for this circular obstacle scenario, while Tables~\ref{tab:free_hyperparameters} and \ref{tab:free_performance} summarize the corresponding experimental setup details and quantitative performance results.

\begin{table}[!ht]
\centering
\resizebox{\textwidth}{!}{%
\begin{tabular}{c|ccccc}
\hline
\textbf{No. of drones} & Drone radius & Perturbation radius & TRS Solver & Annealing & Weight decay $(\lambda)$ \\
\hline
4  & 0.020  & 0.10 & LQR& No & No   \\
8  & 0.020  & 0.10  & LQR + Pretrained TSympOCNet & No& No\\
16 & 0.020  & 0.05 & LQR + Pretrained TSympOCNet & No& No  \\
32 & 0.020  & 0.025 & LQR + Pretrained TSympOCNet& No& Yes $(1e-5)$\\
56 & 0.016  & 0.025 & LQR + Pretrained TSympOCNet& Yes& Yes $(1e-6)$\\
64 & 0.016  & 0.025 & LQR + Pretrained TSympOCNet& Yes& Yes $(1e-6)$ \\
\hline
\end{tabular}
}
\caption{\textbf{Hyperparameters and physical settings for scalability experiments in the presence of a circular obstacle.} The table details the configuration for swarms of varying sizes $N$ navigating around a static circular obstacle. As the agent density increases ($N \ge 32$), the drone radius and initial perturbation radius are reduced to maintain geometric feasibility. Additionally, training stability is enforced through annealing, weight decay ($\lambda$).}
\label{tab:free_hyperparameters}
\end{table}

\begin{table}[!ht]
\resizebox{\textwidth}{!}{%
\centering
\begin{tabular}{c c c c c c c}
\toprule
\# Drones & Pass (Train) & Pass (Test) & Avg PMP (Train) & Avg PMP (Test) & Avg $L$ (Train) & Avg $L$ (Test) \\
\midrule
4  & 100/100 & 100/100 & $3.15 \times 10^{-6}$ & $3.42 \times 10^{-6}$ & 0.579 & 0.580 \\
8  & 100/100 & 94/100  & $4.82 \times 10^{-6}$ & $1.29 \times 10^{-1}$ & 0.886 & 0.890 \\
16 & 100/100 & 99/100  & $2.47 \times 10^{-6}$ & $6.18 \times 10^{-3}$ & 1.625 & 1.623 \\
32 & 100/100 & 97/100  & $3.91 \times 10^{-4}$ & $7.17 \times 10^{-4}$ & 3.806 & 3.807 \\
56 & 100/100 & 99/100  & $3.60 \times 10^{-4}$ & $3.83 \times 10^{-4}$ & 6.275 & 6.276 \\
64 & 100/100 & 98/100  & $7.15 \times 10^{-4}$ & $7.60 \times 10^{-4}$ & 8.770  & 8.762 \\
\bottomrule
\end{tabular}}
\caption{\textbf{Performance of PI-SONet for the `obstacle' case with perturbed initial position.} PI-SONet maintains high feasibility and optimality across all tested swarm sizes.}
\label{tab:free_performance}
\end{table}

\subsubsection{Crossing a square with a maze: handling nonconvex constraints}

We evaluate PI-SONet on a challenging nonconvex optimal control problem -- multi-agent navigation through a maze. Unlike the convex obstacle setting, the maze introduces complex nonconvex and nonsmooth topological constraints that require agents to traverse narrow corridors and avoid dead-ends while maintaining collision-free trajectories.

We consider swarms of $N\in\{4,8\}$ agents navigating from perturbed initial positions to a fixed target configuration. For both swarm sizes, perturbations are sampled uniformly from a ball of radius $0.05$ around a nominal initial configuration. We generate 100 training and 100 testing samples.

Due to the highly nonconvex environment geometry, LQR-based latent trajectories alone are insufficient to capture the correct solution topology. To address this issue, we construct reference trajectories using an Eikonal solver for the unperturbed initial configuration. We observe that the resulting trajectories do indeed encode the correct homotopy class (i.e., feasible paths through the maze). 

\begin{figure}[!ht]
    \centering
    \includegraphics[width=\linewidth]{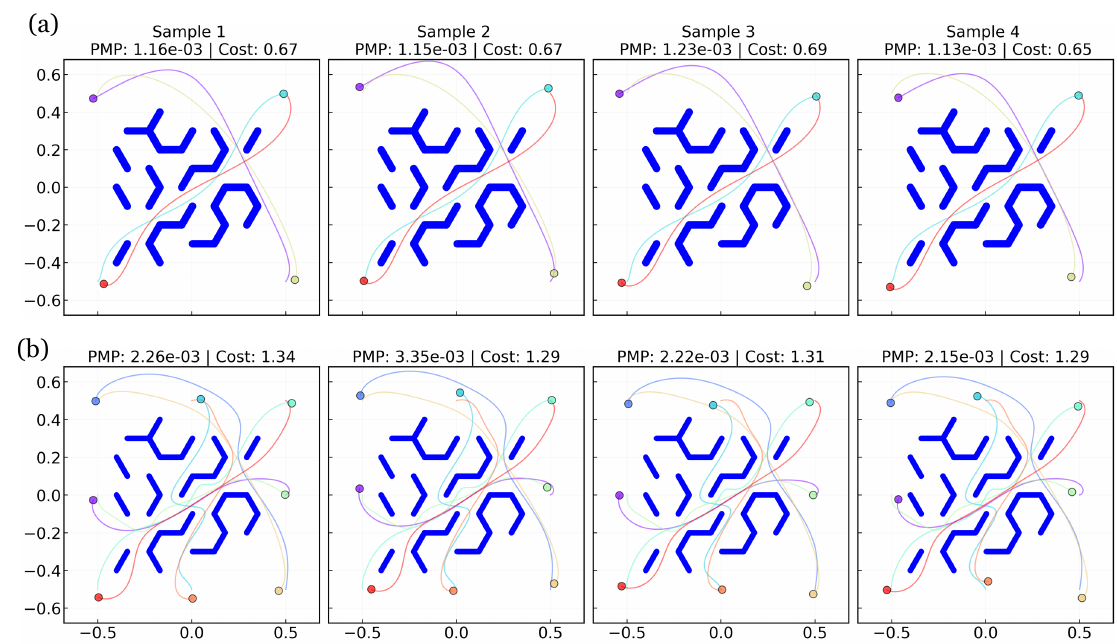}
    \caption{\textbf{Multi-agent navigation in a non-convex maze.} \textbf{(a)} Sample trajectories for 4 agents under perturbed initial conditions. The agents successfully navigate the \textcolor{blue}{blue} obstacles (maze walls) to reach the target configuration. \textbf{(b)} Results for 8-agent swarm under perturbed initial conditions. Despite the high density and narrow corridors, PI-SONet, aided by adaptive activation and annealing, generates collision-free trajectories for the majority of test cases. 
    }\label{fig:sonet_maze}
\end{figure}

We first pretrain the SONet to regress these reference trajectories (positions only, ignoring velocity and co-state) using 500 Adam steps. This stage provides a good initialization that ensures that the predicted trajectories lie in the correct topological class, although they may be suboptimal and may still contain collisions.

Starting from this initialization, we further train the model using the PMP loss to enforce optimality and collision avoidance. In both cases, we do not employ annealing. Instead, we use fixed penalty parameters: $\epsilon = \lambda = 10^{-4}$ for 4 drones and $\epsilon = \lambda = 10^{-3}$ for 8 drones.

The SONet architecture in this case consists of 3 up layers and 3 down layers. The network width for $K$ and $b$ is set to 4 for 4 drones and 8 for 8 drones. For the 4-agent case, the model is trained using 200 Adam steps followed by 200 L-BFGS steps. For the 8-agent case, we use 300 Adam steps followed by 100 L-BFGS steps.

Table~\ref{tab:maze} summarizes the qualitative performance in this maze setting. We observe that PI-SONet achieves strong performance for both swarm sizes with high collision-free success rates on both the training and test sets. For 4 agents, the model achieves 100\% success on the training set and 97\% on the test set. For 8 agents, the model achieves 100\% training success and 95\% test success.

Inference remains efficient even in this non-convex setting. On an NVIDIA H100 GPU, the average batched forward-pass time is $0.022$s for $N=4$ and $0.020$s for $N=8$. These results indicate that the added topological complexity of the maze does not significantly impact runtime, allowing fast trajectory generation even in environments with narrow corridors and constrained feasible regions.

Moreover, the results in this setting highlight the importance of incorporating topology-aware latent trajectories in environments with complex geometries. In particular, the Eikonal-based pretraining effectively guides the model toward feasible homotopy classes, while the PMP-based refinement enforces optimality and collision avoidance. Despite the increased difficulty due to narrow corridors and dense agent interactions, PI-SONet successfully learns a robust solution operator for this highly constrained setting.

\begin{table}[!ht]
\resizebox{\textwidth}{!}{%
\centering
\begin{tabular}{c c c c c c c}
\toprule
\# Drones & Pass (Train) & Pass (Test) & Avg PMP (Train) & Avg PMP (Test) & Avg $L$ (Train) & Avg $L$ (Test) \\
\midrule
4  & 100/100 & 97/100 & $1.17 \times 10^{-3}$ & $ 1.04 \times 10^{-2}$ & 0.667 & 0.668 \\
8  & 100/100 & 95/50 & $2.23 \times 10^{-3}$ & $2.60 \times 10^{-3}$ & 1.295 & 1.301 \\
\bottomrule
\end{tabular}}
\caption{\textbf{Performance of PI-SONet for swarms navigating a 2D maze with perturbed initial positions.} Despite the complex environment geometry, PI-SONet exhibits high feasibility and optimality, demonstrating the power of topology-aware latent right-space solvers and physics-informed training losses.}
\label{tab:maze}
\end{table}

\subsection{Perturbed obstacle geometries}
We evaluate the ability of PI-SONet to approximate the solution operator across varying geometric constraints. Specifically, we consider a multi-agent optimal control problem in which $N\in\{4,8,16\}$ drones navigate around a central circular obstacle of variable radius $r$.

To isolate the effect of the obstacle geometry, the initial and terminal positions of the agents are fixed across all samples, while the obstacle radius is sampled as $r \sim \mathcal{U}[0.05, 0.25]$. For each setting, we generate 50 training and 50 testing samples, where each sample corresponds to a different obstacle size. The latent trajectories are generated using LQR with rotation with parameters $C_B = -\pi/20$ and $C_Q = 1$.

All models use 3 up layers and 3 down layers. The network width for $K$ and $b$ is set to 8 for the 4- and 8-drone cases and 16 for the 16-drone case. Each drone has a fixed radius $C_r = 0.02$. No explicit regularization is used during training. Models are trained for 100 epochs using Adam, followed by 200 epochs using L-BFGS. We employ an annealing strategy for the penalty parameters $\epsilon$ and $\lambda$, initialized at 0.1 and decayed by a factor of 0.6 every 20 epochs.

\begin{table}[!ht]
\resizebox{\textwidth}{!}{%
\centering
\begin{tabular}{c c c c c c c}
\toprule
\# Drones & Pass (Train) & Pass (Test) & Avg PMP (Train) & Avg PMP (Test) & Avg $L$ (Train) & Avg $L$ (Test) \\
\midrule
4  & 50/50 & 50/50 & $1.24 \times 10^{-4}$ & $2.30 \times 10^{-4}$ & 0.556 & 0.559 \\
8  & 50/50 & 50/50 & $4.53 \times 10^{-4}$ & $5.40 \times 10^{-4}$ & 0.838 & 0.846 \\
16 & 50/50 & 48/50 & $4.65 \times 10^{-4}$ & $5.38 \times 10^{-4}$ & 1.549 & 1.565 \\
\bottomrule
\end{tabular}}
\caption{\textbf{Performance of PI-SONet for swarms navigating around a circular obstacle of perturbed radius.} PI-SONet achieves high optimality and feasibility, indicating its ability to capture abstract problem features, such as obstacle geometry.}
\label{tab:heterogeneous_radii}
\end{table}

Table~\ref{tab:heterogeneous_radii} summarizes the performance across different agent counts. PI-SONet achieves perfect or near-perfect feasibility on both training and test sets for 4 and 8 drones with only a slight drop for 16 drones (48/50). The PMP residuals remain low across all cases, indicating accurate satisfaction of optimality conditions. As expected, the average cost increases with the number of agents due to higher interaction complexity and tighter spatial constraints.

Inference remains fast across varying obstacle geometries. On an NVIDIA H100 GPU, the average batched forward-pass time is $0.019$s for $N=4$, $0.019$s for $N=8$, and $0.019$s for $N=16$. These results indicate that changes in obstacle size do not introduce any noticeable computational overhead, allowing rapid trajectory generation even as the feasible region undergoes significant geometric variation.

\begin{figure}
    \centering
    \includegraphics[width=\linewidth]{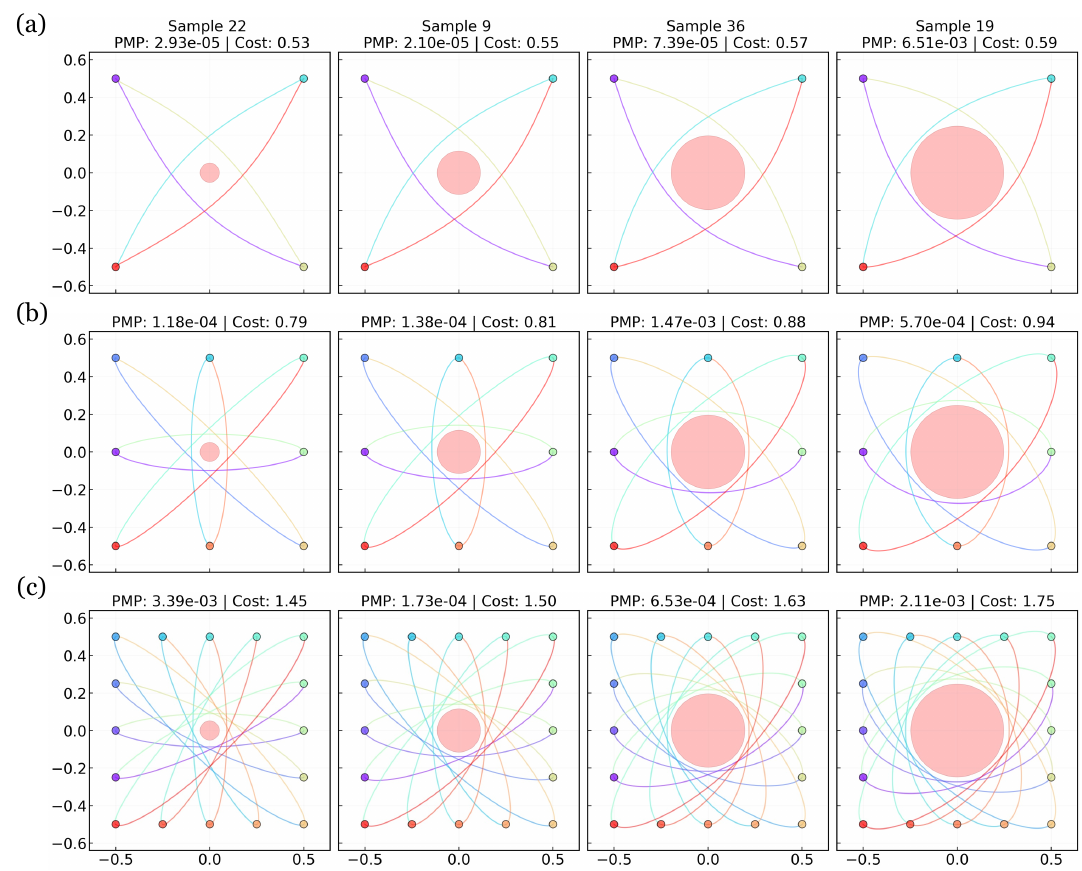}
    \caption{\textbf{PI-SONet under varying obstacle geometries.} We consider drones navigating around a central obstacle with radius $r \in [0.05, 0.25]$. Each column corresponds to a different obstacle radius (increasing from left to right), while rows (a)–(c) denote increasing number of drones, 4, 8 and 16. The initial positions are fixed across all cases to isolate the effect of geometry. PI-SONet successfully adapts to changing obstacle sizes, maintaining feasible, collision-free trajectories across most configurations, with increasing cost as the obstacle becomes larger and constrains the feasible space.}
    \label{fig:heterogenous_radii}
\end{figure}

Figure~\ref{fig:heterogenous_radii} illustrates representative trajectories under varying obstacle sizes. PI-SONet adapts effectively to changing geometries, producing collision-free trajectories across most configurations. As the obstacle radius increases, the feasible region shrinks, leading to more curved trajectories and higher costs, while still preserving feasibility in the majority of cases.

\subsection{Perturbed drone radius: heterogenous setting}
To test the ability of PI-SONet to handle settings with heterogeneous capabilities, we consider scenarios where each drone in the swarm has a different radius.

\subsubsection{Two-dimensional room crossing}

We consider a two-dimensional room-crossing task with heterogeneous drone radii. The dynamics for each drone is now given by
\begin{equation}
\dot{\bx}_i = \bv_i, \quad \dot{\bv}_i = \bu_i - k_i \bv_i \|\bv_i\|_2,
\end{equation}
where the resistance coefficient $k_i$ depends on the drone radius $r_i$. Using the drag relation $F_d \propto \rho C_d A v^2$ and noting that in the 2D setting the projected area scales as $A \propto r_i$ (top-down footprint) while the mass scales as $m_i \propto r_i^2$, we obtain
\[
m_i k_i \propto A \;\;\Rightarrow\;\; r_i^2 k_i \propto r_i \;\;\Rightarrow\;\; k_i \propto \frac{1}{r_i}.
\]
Intuitively, larger drones have more inertia (mass grows as $r_i^2$) relative to the drag they experience (which grows only linearly with $r_i$), leading to a smaller effective resistance coefficient. Conversely, smaller drones experience stronger resistance coefficient, resulting in more velocity-aligned controls.

To evaluate generalization across heterogeneous capabilities, we sample drone radii as $r_i \sim \mathcal{U}[0.01, 0.10]$ for 4 drones and $r_i \sim \mathcal{U}[0.01, 0.05]$ for 8 and 16 drones. For each setting, we generate 50 training and 50 testing samples, where each sample corresponds to a different combination of radii. The obstacle radius is fixed at $C_o = 0.15$ in all cases.

We use LQR with rotation as the latent trajectory solver with parameters $C_B = \pi/20$ and $C_Q = 1.0$. The SONet architecture consists of 3 up layers and 3 down layers with a network width of 64 for both $K$ and $b$. No explicit regularization is used. The resistance coefficient is implemented as $k_i = 1/(50 r_i)$.

Models are trained using a combination of Adam and L-BFGS. For 4 drones, we use 100 Adam epochs followed by 100 L-BFGS epochs, for 8 drones, 100 Adam and 200 L-BFGS epochs, and for 16 drones, 100 Adam and 250 L-BFGS epochs. We employ an annealing strategy for the penalty parameters $(\epsilon, l)$. For 4 and 8 drones, both are initialized at 0.1 and decayed by a factor of 0.6 every 20 epochs. For 16 drones, they are initialized at 0.01 and decayed by a factor of 0.8 every 20 epochs.

Table~\ref{tab:heterogeneous_drone_radii} summarizes the qualitative performance in this scenario across different agent counts. PI-SONet achieves near-perfect feasibility for 4 and 8 drones with only slight degradation for 16 drones due to the corresponding increased agent-interaction complexity. The PMP residuals remain low across all cases, indicating accurate satisfaction of optimality conditions. 

Inference remains efficient under heterogeneous agent dynamics. On an NVIDIA H100 GPU, the average batched forward-pass time is $0.016$s for $N=4$, $0.018$s for $N=8$, and $0.017$s for $N=16$. These results indicate that variations in agent-specific parameters (e.g., here, radius-dependent resistance coefficients) do not significantly impact runtime, indicating that a pre-trained PI-SONet could be trivially transferred in real-time between swarms of slightly different capabilities.

\begin{figure}[!ht]
    \centering
    \includegraphics[width=\linewidth]{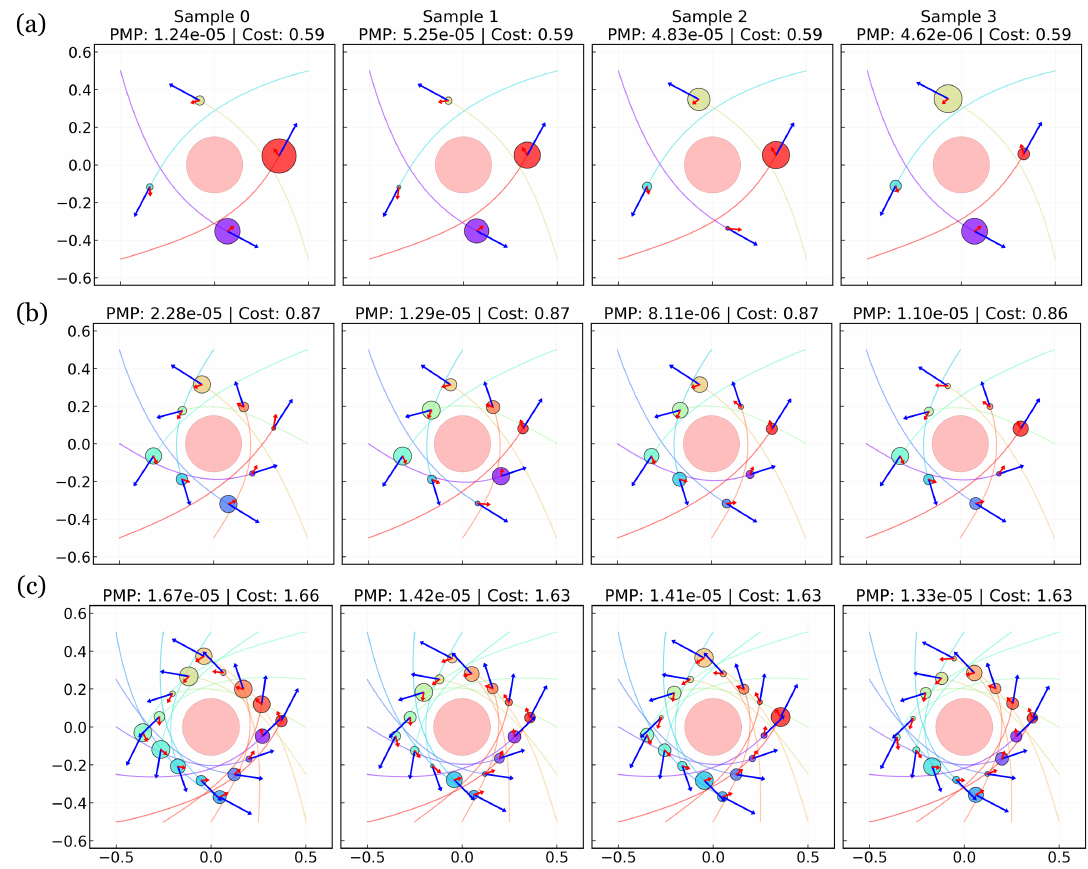}
    \caption{\textbf{PI-SONet for swarms with heterogeneous drone radii.} Snapshots at fixed time $t=6.7$ for (a) 4, (b) 8, and (c) 16 drones. Each column corresponds to a different random realization of radii $r_i \sim \mathcal{U}[0.01, 0.10]$ for 4 drones and $r_i \sim \mathcal{U}[0.01, 0.05]$ for 8 and 16 drones. The obstacle is shown in \textcolor{red}{red}, while colored curves denote agent trajectories. PI-SONet adapts to heterogeneous agent sizes, producing collision-free trajectories in most cases. As the radius decreases (larger $k_i$), the control becomes more aligned with the velocity, reflecting the expected physical response to larger resistance coefficients.}
    \label{fig:varying_drone_radius_supp}
\end{figure}

Figure~\ref{fig:varying_drone_radius_supp} shows representative trajectories at a fixed time. PI-SONet adapts effectively to heterogeneous radii, producing collision-free trajectories in most cases. As the radius decreases (larger $k_i$), the control becomes increasingly aligned with the velocity, reflecting stronger damping effects.

\begin{table}[t]
\resizebox{\textwidth}{!}{%
\centering
\begin{tabular}{c c c c c c c}
\toprule
\# Drones & Pass (Train) & Pass (Test) & Avg PMP (Train) & Avg PMP (Test) & Avg $L$ (Train) & Avg $L$ (Test) \\
\midrule
4  & 50/50 & 49/50 & $1.15 \times 10^{-5}$ & $ 1.71 \times 10^{-5}$ & 0.591 & 0.591 \\
8  & 50/50 & 50/50 & $1.06 \times 10^{-5}$ & $2.85 \times 10^{-5}$ & 0.865 & 0.865 \\
16 & 50/50 & 48/50 & $1.18 \times 10^{-5}$ & $3.49 \times 10^{-2}$ & 1.640 & 1.640 \\
\bottomrule
\end{tabular}}
\caption{\textbf{Performance of PI-SONet for swarms of heterogeneous perturbed drone radii navigating around a circular obstacle in 2D.} PI-SONet demonstrates high feasibility and high optimality, indicating its ability to generalize to heterogeneous agent capabilities.}
\label{tab:heterogeneous_drone_radii}
\end{table}

\subsubsection{Three-dimensional swarm}

We consider a three-dimensional swarm navigation task with $N=100$ drones in the presence of two rectangular obstacles to demonstrate that PI-SONet can be applied to very high-dimensional cases that move beyond planar physical spaces. In this scenario, the agents must move from their initial positions to their target locations while avoiding collisions with both rectangular obstacles in the center of the domain and other drones. The obstacles are defined as axis-aligned boxes:
\begin{align}
\text{Obstacle 1:} \quad & -2 \le x \le 2,\; -0.5 \le y \le 0.5,\; 0 \le z \le 7, \\
\text{Obstacle 2:} \quad & 2 \le x \le 4,\; -1 \le y \le 1,\; 0 \le z \le 4.
\end{align}

The dynamics are linear:
\begin{equation}
\dot{\bx}_i = \bv_i, \quad \dot{\bv}_i = \bu_i,
\end{equation}
and the objective is to minimize the control effort $\int \tfrac{1}{2} \|\bu_i\|^2 \, dt$.

To evaluate robustness to heterogeneous agent sizes, the drone radii are sampled as $r_i \sim \mathcal{U}[0.1, 0.2]$. We generate 50 training and 50 testing samples, where each sample corresponds to a different realization of drone radii. 

The SONet architecture consists of 3 up layers and 3 down layers, with a hidden dimension of 32 for the networks corresponding to $K$ and $b$. We use LQR-generated latent trajectories. To improve training stability, we first train a SONet model for the homogeneous-radius setting and use its weights to initialize the model for the heterogeneous case. This warm-start significantly stabilizes optimization in the high-dimensional, dense-interaction regime.

The model is trained using 100 Adam epochs followed by 50 L-BFGS epochs. We employ annealing for the penalty parameters with $\epsilon$ initialized at $5 \times 10^{-2}$ and $l$ at $5 \times 10^{-3}$, both decayed by a factor of 0.8 every 100 epochs.

PI-SONet successfully scales to large swarms with heterogeneous agent sizes, generating collision-free trajectories while maintaining low control costs. The results demonstrate that our approach mitigates the curse of dimensionality and that the resulting operator can generalize across both high-dimensional state spaces and varying physical agent properties.

\section{Additional benchmarking details}

\subsection{State-of-the-art baseline methods}
For benchmarking, we compare the performance of PI-SONet against the following representative state-of-the-art baseline methods:
\begin{enumerate}
    \item \textit{Multiple shooting method}, which discretizes the optimal control problem in time and solves it as a series of initial-value problems; multiple shooting generally has relatively fast computational times but is highly sensitive to initialization and may struggle with long time horizons and high dimensions.
    \item \textit{Pseudospectral optimal control}, which employes collocation/quadrature to provide high-accuracy solutions but suffers from the curse of dimensionality and generally cannot provide real-time computations.
    \item \textit{TSympOCNet}, which is a sympletic PINN shown to be able scale to very high dimensions/swarm sizes but only learns single-instance solutions and hence does not provide a solution that generalizes to ranges of problem parameters.
\end{enumerate}

For multiple shooting, we use the Python implementation provided by CasADi~\cite{SWcasadi,andersson2019casadi}. We initialize multiple shooting with a linear interpolation between the initial and terminal states, where we add Gaussian noise (mean 0 and standard deviation 0.01) to the velocity to encourage it not to remain at 0 for the whole trajectory. The time to compute the initial guess is not included in the timing results. We set the max number of iterations to be 200 but otherwise use the default CasADi parameters.

For pseudospectral optimal control, we use the MATLAB implementation provided by GPOPS-II~\cite{SWgpops,patterson2014gpops}. We initialize the method with just the initial and terminal states. We set the max number of IPOPT iterations to be 500, the mesh tolerance to be 1e-3, and an initial mesh consisting of 10 sub-meshes with 6 nodes each. 

For TSympOCNet~\cite{zhang2025time}, we use the same JAX-enabled Python implementation that was used in the original paper. We initialize the method using the solution to an SDE system that represents the swarm as a particle system. The time to compute this initialization is not included in the timing results.

A comparison of features between all methods considered is summarized in Table~\ref{tab:benchmark-feature-comparison}. We also list our key innovations over TSympOCNet in Table~\ref{tab:benchmark-pisonet-features}.


\begin{table}[!ht]
\resizebox{\textwidth}{!}{%
\centering
\begin{tabular}{lcccc}
\toprule
 & \textbf{Multiple Shooting} & \textbf{Pseudospectral Method} & \textbf{TSympOCNet} & \textbf{PI-SONet} \\
\midrule
\textit{\makecell[r]{Structure-preserving}} & \xmark & \xmark & \cmark &  \cmark \\ 
\midrule
\textit{\makecell[r]{Continuous constraint \\enforcement}} & \xmark & \xmark& \cmark & \cmark \\
\midrule
\textit{\makecell[r]{Real-time performance
}}& \xmark & \xmark & \xmark & \cmark \\
\midrule
\textit{\makecell[r]{Generalizabilty across \\ problem instances}}& \xmark & \xmark & \xmark & \cmark\\
\bottomrule
\end{tabular}}
\caption{\textbf{Feature comparison between PI-SONet and representative baseline methods.} Like TSympOCNet, PI-SONet preserves symplecticity to maintain optimality and continuously enforces constraints to maintain feasibility. Unlike all baseline methods considered, PI-SONet's generalizability to new problem instances circumvents the need for complete rerunning/retraining, thereby enabling sub-second level inference times  on new problem settings.}
\label{tab:benchmark-feature-comparison}
\end{table}

\begin{table}[!ht]
\centering
\begin{tabular}{ll}
\toprule
\textbf{TSympOCNet} & \textbf{PI-SONet (ours)} \\
\midrule

\cmark\ Parameterized symplectic map
&
\cmark\ Parameterized \textit{family} of symplectic maps \\[0.5em]

\cmark \makecell[l]{Families of TRS solvers: \\ LQR, LQR + TL-SympNet}
&
\cmark\makecell[l]{New families of TRS solvers:\\ LQR, LQR with rotation, LQR + pretrained TSympOCNet, Eikonal} \\[1.0em]

\text{\sffamily \textcolor{red}{X}} Single-instance solver &
\cmark\ Generalizable to new family of problem scenarios \\[0.5em]

\text{\sffamily \textcolor{red}{X}} Retrain for new problem scenarios &
\cmark\ Fast test-time finetuning without re-training \\[0.5em]

\text{\sffamily \textcolor{red}{X}} Full $K$ matrices &
\cmark\ Efficient block diagonal structure for $K$ matrix \\

\bottomrule
\end{tabular}
\caption{\textbf{Comparison between TSympOCNet and PI-SONet.}
While TSympOCNet learns a parameterized symplectic map for a fixed problem instance, PI-SONet extends this framework to a parameterized \textit{family} of symplectic maps, enabling generalization across classes of problem scenarios. PI-SONet performs real-time inference without retraining and fast test-time fine tuning when needed, incorporates richer families of latent TRS solvers, and leverages an efficient block-diagonal structure for the learned operators for improved scalability.}
\label{tab:benchmark-pisonet-features}
\end{table}
\subsection{Some qualitative results for the \textit{Maze} scenario}

\begin{figure}[!ht]
    \centering
    \begin{subfigure}[b]{0.4\textwidth}
        \includegraphics[width=\textwidth]{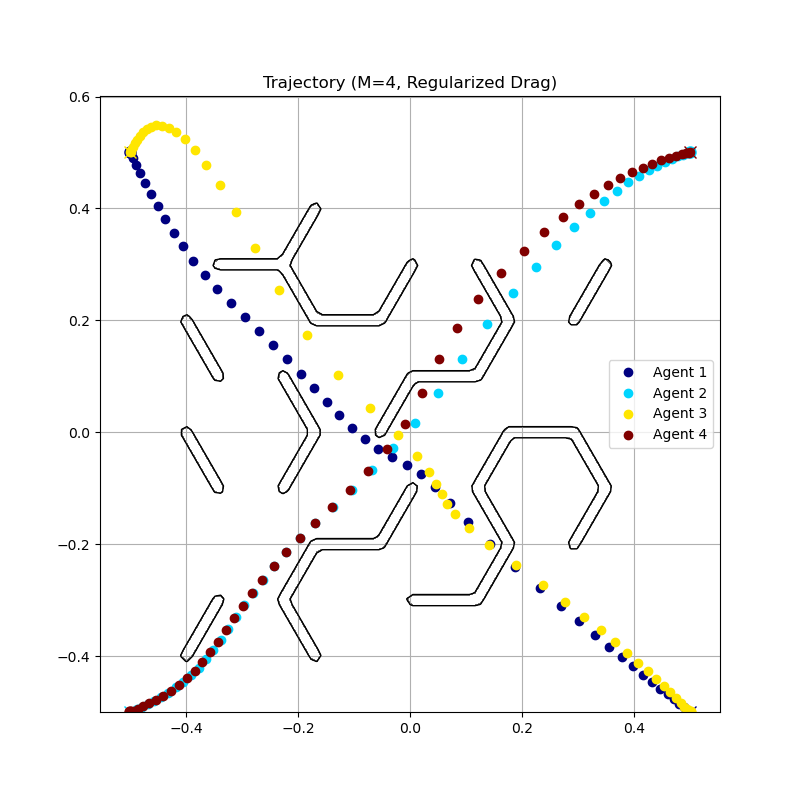}
        \caption{Discretized Trajectory}
    \end{subfigure}
    \begin{subfigure}[b]{0.4\textwidth}
        \includegraphics[width=\textwidth]{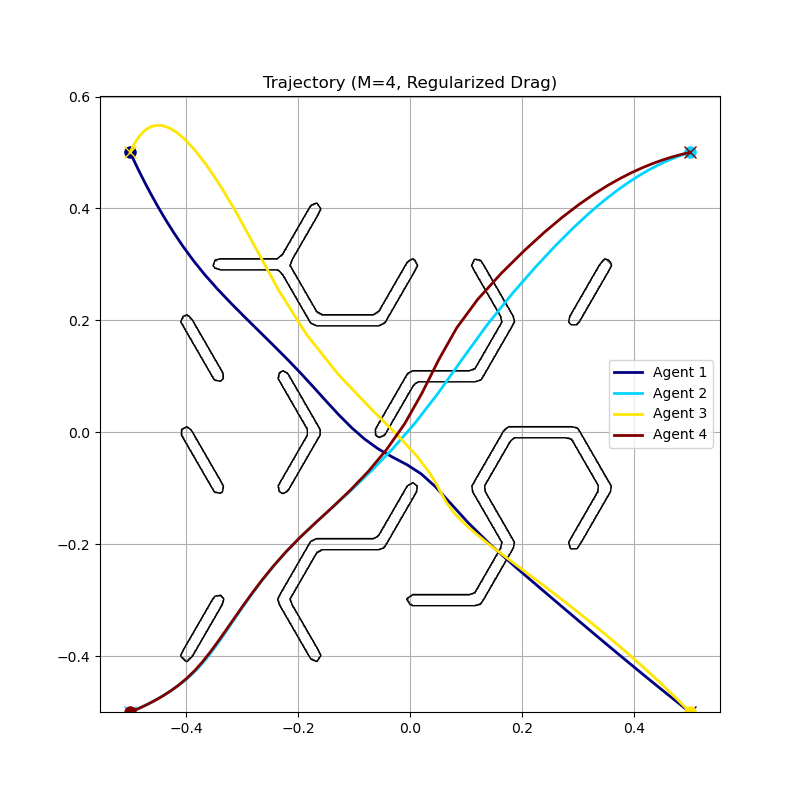}
        \caption{Interpolated Trajectory}
    \end{subfigure}
    \caption{\textbf{Sample trajectories using discretization-based baseline methods for a 4-agent swarm navigating through a nonconvex maze.} Without a good initialization, discretization-based methods, such as multiple shooting and pseudospectral method, often lead to infeasible trajectories as problem constraints are only enforced at discrete points. This behavior exemplifies the issue of initialization sensitivity that many existing baseline methods experience and which becomes more pronounced in nonconvex cases, such as the one shown here. In such cases, using a finer grid may not be enough to recover convergence to a feasible solution (and in fact, it was not able to in the case shown here). In contrast, as an operator-based approach, PI-SONet continuously enforces constraints for all times and states and generally does not suffer from this behavior.}
    \label{fig:benchmark_maze}
\end{figure}

Discretization-based methods, such as multiple shooting and pseudospectral method, generally rely on having a good initial guess in order to converge to a feasible, (sub)optimal solution. However, problem constraints (e.g., collision or obstacle avoidance constraints) are only enforced at discrete grid points. Hence, when not initialized properly, these methods may report generated trajectories as successful when, in actuality, they contain visually obvious collisions. An example of this behavior is shown in Figure~\ref{fig:benchmark_maze}, in which multiple shooting generates trajectories with no collisions at the discretization points, but any reasonable interpolation of this discrete trajectory results in obvious and severe crashes with the maze.

\section{Ablation studies}

\subsection{Choice of latent solver}
\label{subsec:latent_solver}
We investigate the effect of the latent solver used to generate the reference trajectories. In all cases, the downstream PI-SONet architecture and training procedure are kept fixed, and only the latent trajectories are varied. We consider three choices: (i) standard LQR, (ii) LQR with an additional rotation term in the latent dynamics, and (iii) LQR composed with a pretrained TSympOCNet. The pretrained TSympOCNet is trained only on the unperturbed initial condition and is then reused to generate latent trajectories for perturbed cases. All experiments in this section are conducted in the \emph{Free} setting, i.e., without obstacles.

\begin{figure}
    \centering
    \includegraphics[width=\linewidth]{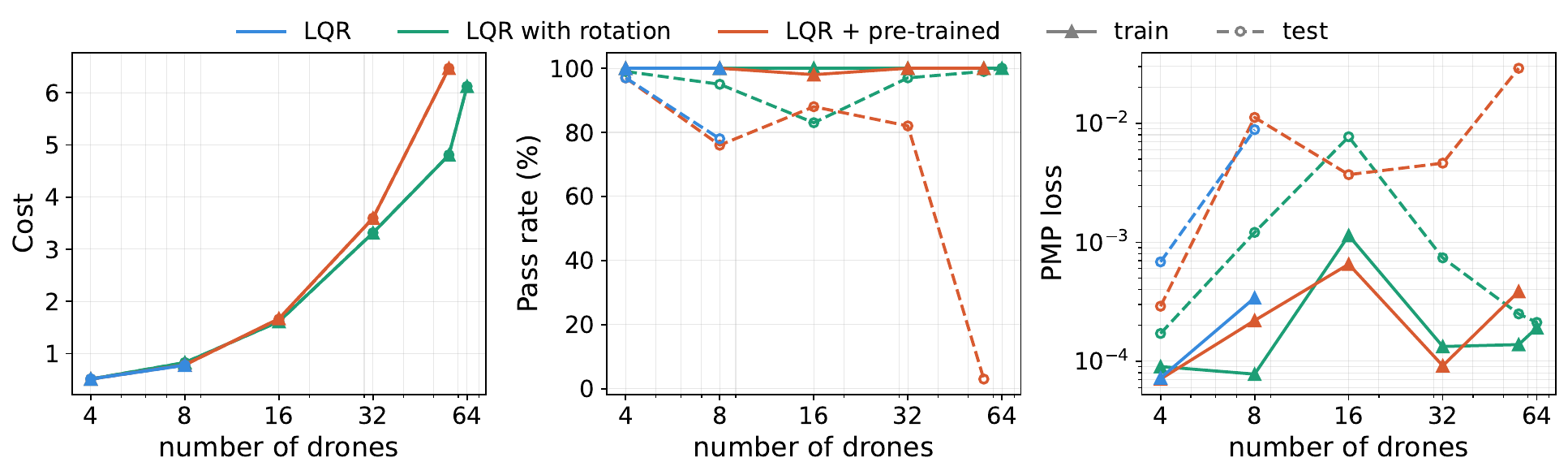}
    \caption{\textbf{Effect of different latent solvers (in the `free' case).}
    We compare LQR, LQR with rotation, and LQR composed with a pretrained TSympOCNet across increasing swarm sizes in the absence of obstacles. 
    (\textit{Left}) Cost increases with the number of drones for all methods with LQR with rotation consistently achieving lower costs than LQR + pretrained TSympOCNet at larger scales. 
    (\textit{Middle}) Pass rate (\%) on train (solid) and test (dashed) datasets. While all methods perform similarly for small systems, LQR with rotation maintains high success rates as the number of drones increases, whereas LQR + pretrained TSympOCNet exhibits a sharp degradation at larger scales. 
    (\textit{Right}) PMP residuals (log scale). LQR with rotation generally achieves lower test residuals compared to the other methods, indicating better adherence to optimality conditions.}
    \label{fig:placeholder}
\end{figure}

The quantitative results are summarized in Table~\ref{tab:drone_scaling} with corresponding trends shown in Figure~\ref{fig:placeholder}. For small systems ($n=4,8$), all three approaches achieve comparable performance in terms of cost and feasibility. However, as the number of drones increases, clear differences emerge. The standard LQR solver fails to train stably beyond small system sizes, indicating that the latent trajectories do not provide a suitable initialization for the downstream optimization.

In contrast, incorporating rotation into the latent dynamics significantly improves scalability. LQR with rotation maintains near-perfect pass rates across all tested system sizes with consistently low PMP residuals and competitive objective costs. These results suggest that even a simple modification to the latent dynamics can better capture the structure of the solution manifold, leading to improved robustness as the system dimension grows.

The pretrained TSympOCNet composition performs competitively at small scales, but its performance deteriorates as the number of drones increases. In particular, we observe a drop in pass rate and a sharp increase in PMP residuals for larger systems (e.g., $n=56$), indicating poor generalization. This behavior is expected, as the pretrained model is trained only on the unperturbed configuration and is not adapted to variations in initial conditions or system size.

Overall, these results highlight the importance of the choice of latent solver in enabling scalable learning of optimal control solutions. Even in this geometrically simple free-space setting, incorporating rotation in the latent dynamics significantly improves robustness and scalability, while naive reuse of a pretrained model fails to generalize to larger systems.

\begin{table}[t]
\resizebox{\textwidth}{!}{
\centering
\small
\setlength{\tabcolsep}{5pt}
\begin{tabular}{llcccccc}
\toprule
 & & \multicolumn{3}{c}{Train} & \multicolumn{3}{c}{Test} \\
\cmidrule(lr){3-5} \cmidrule(lr){6-8}
Metric & $n$ & LQR & LQR w/ rotation & LQR + pretrained & LQR & LQR w/ rotation & LQR + pretrained \\
\midrule
\multirow{6}{*}{cost}
  & 4  & 0.512 & 0.511 & 0.511 & 0.513 & 0.512 & 0.511 \\
  & 8  & 0.776 & 0.826 & 0.782 & 0.778 & 0.828 & 0.783 \\
  & 16 & ---   & 1.614 & 1.670 & ---   & 1.612 & 1.660 \\
  & 32 & ---   & 3.313 & 3.597 & ---   & 3.315 & 3.596 \\
  & 56 & ---   & 4.807 & 6.470 & ---   & 4.809 & 6.470 \\
  & 64 & ---   & 6.123 & ---   & ---   & 6.122 & --- \\
\midrule
\multirow{6}{*}{pass (\%)}
  & 4  & 100 & 100 & 100 & 97 & 99  & 97 \\
  & 8  & 100 & 100 & 100 & 78 & 95  & 76 \\
  & 16 & --- & 100 & 98  & --- & 83 & 88 \\
  & 32 & --- & 100 & 100 & --- & 97 & 82 \\
  & 56 & --- & 100 & 100 & --- & 99 & 3 \\
  & 64 & --- & 100 & --- & --- & 100 & --- \\
\midrule
\multirow{6}{*}{PMP}
  & 4  & $7.24\mathrm{e}{-5}$ & $9.02\mathrm{e}{-5}$ & $7.08\mathrm{e}{-5}$ & $6.84\mathrm{e}{-4}$ & $1.71\mathrm{e}{-4}$ & $2.90\mathrm{e}{-4}$ \\
  & 8  & $3.41\mathrm{e}{-4}$ & $7.79\mathrm{e}{-5}$ & $2.20\mathrm{e}{-4}$ & $8.84\mathrm{e}{-3}$ & $1.21\mathrm{e}{-3}$ & $1.12\mathrm{e}{-2}$ \\
  & 16 & ---                  & $1.14\mathrm{e}{-3}$ & $6.54\mathrm{e}{-4}$ & ---                  & $7.70\mathrm{e}{-3}$ & $3.70\mathrm{e}{-3}$ \\
  & 32 & ---                  & $1.33\mathrm{e}{-4}$ & $9.16\mathrm{e}{-5}$ & ---                  & $7.39\mathrm{e}{-4}$ & $4.62\mathrm{e}{-3}$ \\
  & 56 & ---                  & $1.38\mathrm{e}{-4}$ & $3.86\mathrm{e}{-4}$ & ---                  & $2.50\mathrm{e}{-4}$ & $2.90\mathrm{e}{-2}$ \\
  & 64 & ---                  & $1.91\mathrm{e}{-4}$ & ---                  & ---                  & $2.12\mathrm{e}{-4}$ & --- \\
\bottomrule
\end{tabular}}
\caption{\textbf{Quantitative comparison of latent solvers in free space under increasing system size.}
We report cost, pass rate (\%), and PMP residuals for LQR, LQR with rotation, and LQR composed with a pretrained TSympOCNet, evaluated on both training and test sets without obstacles. 
LQR fails to train stably beyond small system sizes. Incorporating rotation into the latent dynamics significantly improves both feasibility and optimality, maintaining high pass rates and low PMP residuals across all scales. 
In contrast, while the pretrained composition performs competitively at small scales, it exhibits degraded generalization and stability as the number of drones increases, reflected in reduced pass rates and higher PMP residuals.}
\label{tab:drone_scaling}
\end{table}

\subsection{Effect of enforcing symplectic structure}

To assess the importance of enforcing symplectic structure, we replace the SONet architecture with a standard multilayer perceptron (MLP) that directly maps latent trajectories to physical states. Specifically, the MLP takes as inputs the concatenated latent state $(\by(t), \bq(t))$, obstacle radius $r$, and time $t$ and outputs $(\bx'(t), \bp'(t))$. The network consists of fully-connected layers with SiLU activations and a hidden dimension of $d=256$, chosen such that the total number of parameters is comparable to that of PI-SONet across all problem sizes. 

To ensure a fair comparison, we enforce boundary conditions in the same manner as SONet by constructing the output as
\[
\bx_{\text{out}}(t) = \by(t) + t(1-t)\,\bx'(t), \quad 
\bp_{\text{out}}(t) = \bp'(t),
\]
so that the initial and terminal constraints are satisfied by construction. The training procedure, loss formulation, and optimization schedule (Adam followed by L-BFGS) are kept identical to the SONet setting. For the 16-drone case, a small $\ell_2$ regularization ($\lambda = 10^{-6}$) was required to stabilize training.

Table~\ref{tab:symp_ablation} summarizes the quantitative comparison. While both models achieve similar success rates in terms of feasibility, the MLP consistently exhibits higher PMP residuals and suboptimal costs with the gap becoming more pronounced as the number of drones increases. Notably, these differences arise despite comparable parameter counts, indicating that the performance gain of PI-SONet is not simply due to model capacity.

\begin{table*}[!ht]
\resizebox{\textwidth}{!}{
\centering
\small
\begin{tabular}{c|c|cccc|cccc}
\toprule
\textbf{\# Drones} & \textbf{Model} & \textbf{\# Params} & \textbf{Train Time} & \textbf{PMP (Train)} & \textbf{Cost (Train)} & \textbf{Test Time} & \textbf{PMP (Test)} & \textbf{Cost (Test)} & \textbf{Pass} \\
\midrule

\multirow{2}{*}{4}
& MLP   & 82,976  & 25.45  & $2.31{\times}10^{-4}$ & 0.601 & 0.0034 & $1.12{\times}10^{-3}$ & 0.601 & 49/50 \\
& SONet & 76,608  & 106.42 & $1.24{\times}10^{-4}$ & 0.556 & 0.0189 & $7.24{\times}10^{-5}$ & 0.559 & 50/50 \\

\midrule

\multirow{2}{*}{8}
& MLP   & 99,392  & 26.28  & $6.19{\times}10^{-4}$ & 0.967 & 0.0053 & $2.22{\times}10^{-3}$ & 0.968 & 49/50 \\
& SONet & 90,624  & 85.35  & $4.53{\times}10^{-4}$ & 0.838 & 0.0192 & $3.41{\times}10^{-4}$ & 0.846 & 50/50 \\

\midrule

\multirow{2}{*}{16}
& MLP   & 132,224 & 34.26  & $2.10{\times}10^{-4}$ & 1.753 & 0.0046 & $7.28{\times}10^{-4}$ & 1.754 & 49/50 \\
& SONet & 134,208 & 107.74 & $4.65{\times}10^{-4}$ & 1.549 & 0.0191 & $2.30{\times}10^{-4}$ & 1.565 & 48/50 \\

\bottomrule
\end{tabular}}
\caption{\textbf{Ablation study on enforcing symplectic structure}. We compare a standard MLP against SONet (symplectic map). SONet achieves consistently lower cost and PMP residuals with comparable parameter counts, highlighting the importance of preserving symplecticity.}
\label{tab:symp_ablation}
\end{table*}

\begin{figure}[!ht]
    \centering
    \includegraphics[width=\linewidth]{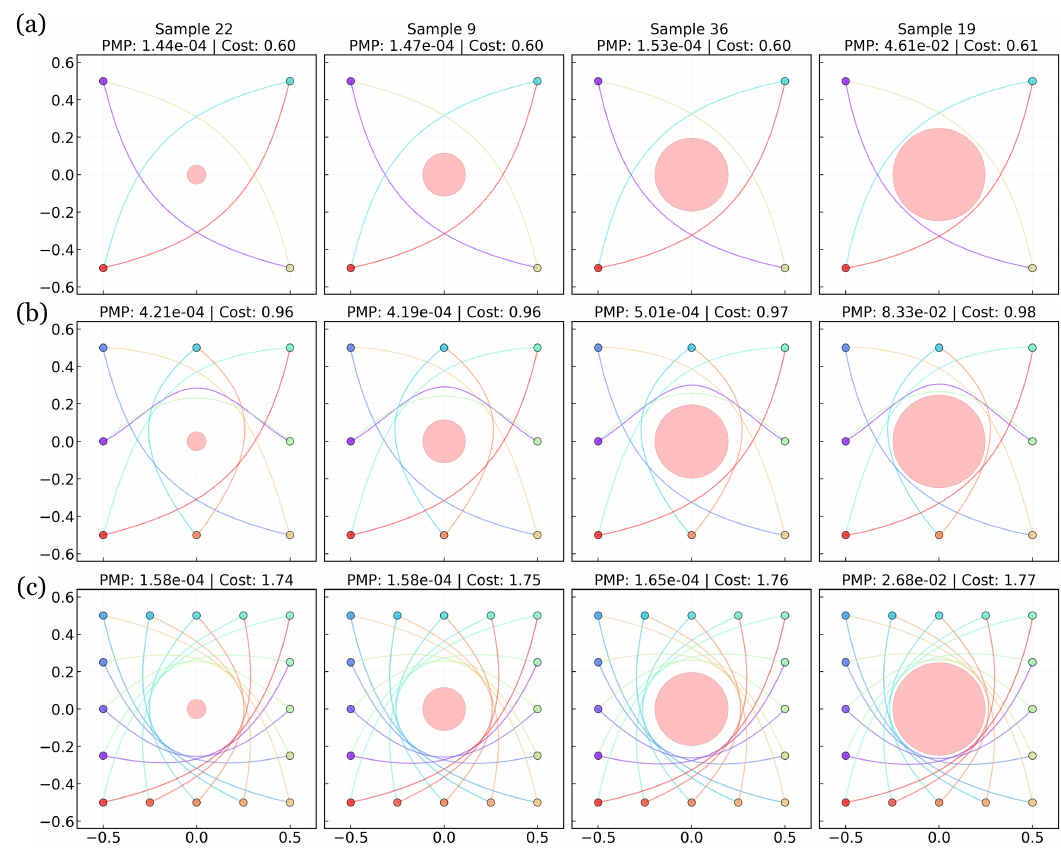}
    \caption{\textbf{Ablation study: MLP trajectories without symplectic structure under varying obstacle geometries.} 
    We consider drones navigating around a central obstacle with radius $r \in [0.05, 0.25]$. Each column corresponds to a different obstacle radius (increasing from left to right), while rows (a)–(c) denote increasing numbers of drones (4, 8 and 16). The initial positions are fixed across all cases to isolate the effect of the obstacle geometry. 
    Unlike PI-SONet, the MLP fails to adapt its trajectories to the changing obstacle size; the predicted paths remain qualitatively similar across columns, leading to increased PMP residuals and costs. These results highlight the importance of enforcing symplectic structure for capturing geometry-dependent optimal solutions.}
    \label{fig:ablation_symplectic}
\end{figure}

The qualitative behavior is illustrated in Figure~\ref{fig:ablation_symplectic}. Unlike PI-SONet, which adapts its trajectories to the obstacle geometry, the MLP produces trajectories that are largely insensitive to changes in the obstacle radius. As the obstacle size increases, the predicted paths remain qualitatively similar, resulting in degraded optimality and increased PMP residuals. Although the trajectories are mostly collision-free, they fail to adjust in a geometry-aware manner, leading to suboptimal paths.

These results highlight that enforcing symplectic structure plays an important role in capturing the underlying Hamiltonian dynamics and geometry-dependence. In its absence, even an analogously expressive MLP struggles to learn the correct dependence on problem parameters, resulting in solutions that satisfy constraints but deviate from optimality.



%% file: ref.bib
@book{bellman1957dynamic,
  author    = {Richard Bellman},
  title     = {Dynamic Programming},
  publisher = {Princeton University Press},
  address   = {Princeton, NJ},
  year      = {1957},
  isbn      = {9780691079516}
}

@article{meng2022sympocnet,
  title={Sympocnet: Solving optimal control problems with applications to high-dimensional multiagent path planning problems},
  author={Meng, Tingwei and Zhang, Zhen and Darbon, Jerome and Karniadakis, George},
  journal={SIAM Journal on Scientific Computing},
  volume={44},
  number={6},
  pages={B1341--B1368},
  year={2022},
  publisher={SIAM}
}

@book{trelat2005controle,
  author    = {Emmanuel Tr{\'e}lat},
  title     = {Contr{\^o}le optimal: th{\'e}orie et applications},
  publisher = {Vuibert},
  year      = {2005},
  series    = {Math{\'e}matiques Concr{\`e}tes},
  address   = {Paris},
  isbn      = {9782711771752}
}

@article{lu2021learning,
  title={Learning nonlinear operators via DeepONet based on the universal approximation theorem of operators},
  author={Lu, Lu and Jin, Pengzhan and Pang, Guofei and Zhang, Zhongqiang and Karniadakis, George Em},
  journal={Nature machine intelligence},
  volume={3},
  number={3},
  pages={218--229},
  year={2021},
  publisher={Nature Publishing Group UK London}
}

@article{darbon2016algorithms,
  title={Algorithms for overcoming the curse of dimensionality for certain Hamilton--Jacobi equations arising in control theory and elsewhere},
  author={Darbon, J{\'e}r{\^o}me and Osher, Stanley},
  journal={Research in the Mathematical Sciences},
  volume={3},
  number={1},
  pages={19},
  year={2016},
  publisher={Springer}
}

@article{darbon2023neural,
  title={Neural network architectures using min-plus algebra for solving certain high-dimensional optimal control problems and Hamilton--Jacobi PDEs},
  author={Darbon, J{\'e}r{\^o}me and Dower, Peter M and Meng, Tingwei},
  journal={Mathematics of Control, Signals, and Systems},
  volume={35},
  number={1},
  pages={1--44},
  year={2023},
  publisher={Springer}
}

@article{wang2021learning,
  title={Learning the solution operator of parametric partial differential equations with physics-informed DeepONets},
  author={Wang, Sifan and Wang, Hanwen and Perdikaris, Paris},
  journal={Science advances},
  volume={7},
  number={40},
  pages={eabi8605},
  year={2021},
  publisher={American Association for the Advancement of Science}
}

@book{lewis2012optimal,
  author    = {Frank L. Lewis and Draguna L. Vrabie and Vassilis L. Syrmos},
  title     = {Optimal Control},
  edition   = {3},
  publisher = {John Wiley \& Sons},
  address   = {Hoboken, NJ},
  year      = {2012},
  isbn      = {9780470633496},
  doi       = {10.1002/9781118122631}
}

@article{li2020fourier,
  title={Fourier neural operator for parametric partial differential equations},
  author={Li, Zongyi and Kovachki, Nikola and Azizzadenesheli, Kamyar and Liu, Burigede and Bhattacharya, Kaushik and Stuart, Andrew and Anandkumar, Anima},
  journal={arXiv preprint arXiv:2010.08895},
  year={2020}
}

@incollection{caillau2023algorithmic,
  author    = {Jean-Baptiste Caillau and Roberto Ferretti and Emmanuel Tr{\'e}lat and Hasnaa Zidani},
  title     = {An algorithmic guide for finite-dimensional optimal control problems},
  booktitle = {Handbook of Numerical Analysis},
  editor    = {Emmanuel Tr{\'e}lat and Enrique Zuazua},
  volume    = {24},
  pages     = {559--626},
  year      = {2023},
  publisher = {North-Holland},
  address   = {Amsterdam},
  doi       = {10.1016/bs.hna.2022.11.006}
}

@article{foderaro2014distributed,
  author  = {Greg Foderaro and Silvia Ferrari and Thomas A. Wettergren},
  title   = {Distributed optimal control for multi-agent trajectory optimization},
  journal = {Automatica},
  volume  = {50},
  number  = {1},
  pages   = {149--154},
  year    = {2014}
}

@article{fahroo2002direct,
  title={Direct trajectory optimization by a Chebyshev pseudospectral method},
  author={Fahroo, Fariba and Ross, I Michael},
  journal={Journal of Guidance, Control, and Dynamics},
  volume={25},
  number={1},
  pages={160--166},
  year={2002}
}

@article{robinson2018efficient,
  title={An efficient algorithm for optimal trajectory generation for heterogeneous multi-agent systems in non-convex environments},
  author={Robinson, D Reed and Mar, Robert T and Estabridis, Katia and Hewer, Gary},
  journal={IEEE Robotics and Automation Letters},
  volume={3},
  number={2},
  pages={1215--1222},
  year={2018},
  publisher={IEEE}
}

@article{gong2006pseudospectral,
  title={A pseudospectral method for the optimal control of constrained feedback linearizable systems},
  author={Gong, Qi and Kang, Wei and Ross, I Michael},
  journal={IEEE transactions on automatic control},
  volume={51},
  number={7},
  pages={1115--1129},
  year={2006},
  publisher={IEEE}
}

@article{boucher2016galerkin,
  title={Galerkin optimal control},
  author={Boucher, Randy and Kang, Wei and Gong, Qi},
  journal={Journal of Optimization Theory and Applications},
  volume={169},
  pages={825--847},
  year={2016},
  publisher={Springer}
}

@inproceedings{kirchner2020hj,
  title={A {Hamilton--Jacobi} Formulation for Optimal Coordination of Heterogeneous Multiple Vehicle Systems},
  author={Kirchner, Matthew R and Debord, Mark J and Hespanha, Jo{\~a}o P},
  booktitle={2020 IEEE/RSJ International Conference on Intelligent Robots and Systems (IROS)},
  pages={11623--11630},
  year={2020},
  organization={IEEE}
}

@article{lee2021hopf,
  title={A {Hopf-Lax} Formula in {Hamilton--Jacobi} Analysis of Reach-Avoid Problems},
  author={Lee, Donggun and Tomlin, Claire J},
  journal={IEEE Control Systems Letters},
  volume={5},
  number={3},
  pages={1055--1060},
  year={2020},
  publisher={IEEE}
}

@Article{zbMATH07547920,
 Author = {Oster, Mathias and Sallandt, Leon and Schneider, Reinhold},
 Title = {Approximating optimal feedback controllers of finite horizon control problems using hierarchical tensor formats},
 FJournal = {SIAM Journal on Scientific Computing},
 Journal = {SIAM J. Sci. Comput.},
 ISSN = {1064-8275},
 Volume = {44},
 Number = {3},
 Pages = {b746--b770},
 Year = {2022},
 Language = {English},
 DOI = {10.1137/21M1412190},
 zbMATH = {7547920}
}

@Article{zbMATH07364328,
 Author = {Dolgov, Sergey and Kalise, Dante and Kunisch, Karl K.},
 Title = {Tensor decomposition methods for high-dimensional {Hamilton}-{Jacobi}-{Bellman} equations},
 FJournal = {SIAM Journal on Scientific Computing},
 Journal = {SIAM J. Sci. Comput.},
 ISSN = {1064-8275},
 Volume = {43},
 Number = {3},
 Pages = {a1625--a1650},
 Year = {2021},
 Language = {English},
 DOI = {10.1137/19M1305136},
 zbMATH = {7364328},
 Zbl = {1471.65184}
}

@article{kunisch2023learning,
  title={Learning optimal feedback operators and their sparse polynomial approximations},
  author={Kunisch, Karl and V{\'a}squez-Varas, Donato and Walter, Daniel},
  journal={Journal of Machine Learning Research},
  volume={24},
  number={301},
  pages={1--38},
  year={2023}
}

@article{Han2016DeepLA,
  title={Deep Learning Approximation for Stochastic Control Problems},
  author={Jiequn Han and E Weinan},
  journal={ArXiv},
  year={2016},
  volume={abs/1611.07422},
  url={https://api.semanticscholar.org/CorpusID:16386889}
}

@article{nakamura2022neural,
  title={Neural network optimal feedback control with guaranteed local stability},
  author={Nakamura-Zimmerer, Tenavi and Gong, Qi and Kang, Wei},
  journal={IEEE Open Journal of Control Systems},
  volume={1},
  pages={210--222},
  year={2022},
  publisher={IEEE}
}

@article{andersson2019casadi,
  title={{CasADi: a software framework for nonlinear optimization and optimal control}},
  author={Andersson, Joel AE and Gillis, Joris and Horn, Greg and Rawlings, James B and Diehl, Moritz},
  journal={Mathematical Programming Computation},
  volume={11},
  number={1},
  pages={1--36},
  year={2019},
  publisher={Springer}
}

@misc{SWcasadi,
  author = {Joel Andersson and Joris Gillis},
  title = {{CasADi}. {V}ersion 3.7.2. \url{https://web.casadi.org/}},
  url = {https://web.casadi.org/},
  version = {3.7.2}
}

@article{patterson2014gpops,
  title={{GPOPS-II: A MATLAB software for solving multiple-phase optimal control problems using hp-adaptive Gaussian quadrature collocation methods and sparse nonlinear programming}},
  author={Patterson, Michael A and Rao, Anil V},
  journal={ACM Transactions on Mathematical Software (TOMS)},
  volume={41},
  number={1},
  pages={1--37},
  year={2014},
  publisher={ACM New York, NY, USA}
}

@misc{SWgpops,
  author = {Michael A. Patterson and Anil V. Rao},
  title = {{GPOPS-II: Next-Generation Optimal Control Software}. {V}ersion 2.3. \url{https://www.gpops2.com/}},
  url = {https://www.gpops2.com/},
  version = {2.3},
}

@inproceedings{mavrogiannis2020decentralized,
  title={Decentralized multi-agent navigation planning with braids},
  author={Mavrogiannis, Christoforos I and Knepper, Ross A},
  booktitle={Algorithmic Foundations of Robotics XII: Proceedings of the Twelfth Workshop on the Algorithmic Foundations of Robotics},
  pages={880--895},
  year={2020},
  organization={Springer}
}

@article{sperl2025separable,
  title={Separable approximations of optimal value functions and their representation by neural networks},
  author={Sperl, Mario and Saluzzi, Luca and Kalise, Dante and Gr{\"u}ne, Lars},
  journal={arXiv preprint arXiv:2502.08559},
  year={2025}
}

@article{akian2026tropical,
  title={Tropical low-rank approximation and application to optimal control of N-body systems},
  author={Akian, Marianne and Gaubert, Stephane and Liu, Shanqing and Qi, Yang},
  journal={arXiv preprint arXiv:2604.18785},
  year={2026}
}

@article{hure2020deep,
  title={Deep backward schemes for high-dimensional nonlinear PDEs},
  author={Hur{\'e}, C{\^o}me and Pham, Huy{\^e}n and Warin, Xavier},
  journal={Mathematics of Computation},
  volume={89},
  number={324},
  pages={1547--1579},
  year={2020}
}

@article{gaspard2025monotone,
  title={Monotone Causality in Opportunistically Stochastic Shortest Path Problems},
  author={Gaspard, Mallory E and Vladimirsky, Alexander},
  journal={Mathematics of Operations Research},
  year={2025},
  publisher={INFORMS}
}

@article{sethian2003ordered,
  title={Ordered upwind methods for static Hamilton--Jacobi equations: Theory and algorithms},
  author={Sethian, James A and Vladimirsky, Alexander},
  journal={SIAM Journal on Numerical Analysis},
  volume={41},
  number={1},
  pages={325--363},
  year={2003},
  publisher={SIAM}
}

@article{dower2015max,
  title={A max-plus dual space fundamental solution for a class of operator differential Riccati equations},
  author={Dower, Peter M and McEneaney, William M},
  journal={SIAM Journal on Control and Optimization},
  volume={53},
  number={2},
  pages={969--1002},
  year={2015},
  publisher={SIAM}
}

@article{brunton2024promising,
  title={Promising directions of machine learning for partial differential equations},
  author={Brunton, Steven L and Kutz, J Nathan},
  journal={Nature Computational Science},
  volume={4},
  number={7},
  pages={483--494},
  year={2024},
  publisher={Nature Publishing Group US New York}
}

@article{yin2024scalable,
  title={A scalable framework for learning the geometry-dependent solution operators of partial differential equations},
  author={Yin, Minglang and Charon, Nicolas and Brody, Ryan and Lu, Lu and Trayanova, Natalia and Maggioni, Mauro},
  journal={Nature computational science},
  volume={4},
  number={12},
  pages={928--940},
  year={2024},
  publisher={Nature Publishing Group US New York}
}

@article{vinuesa2022enhancing,
  title={Enhancing computational fluid dynamics with machine learning},
  author={Vinuesa, Ricardo and Brunton, Steven L},
  journal={Nature Computational Science},
  volume={2},
  number={6},
  pages={358--366},
  year={2022},
  publisher={Nature Publishing Group US New York}
}

@article{raabe2023accelerating,
  title={Accelerating the design of compositionally complex materials via physics-informed artificial intelligence},
  author={Raabe, Dierk and Mianroodi, Jaber Rezaei and Neugebauer, J{\"o}rg},
  journal={Nature computational science},
  volume={3},
  number={3},
  pages={198--209},
  year={2023},
  publisher={Nature Publishing Group US New York}
}

@article{cao2024laplace,
  title={Laplace neural operator for solving differential equations},
  author={Cao, Qianying and Goswami, Somdatta and Karniadakis, George Em},
  journal={Nature Machine Intelligence},
  volume={6},
  number={6},
  pages={631--640},
  year={2024},
  publisher={Nature Publishing Group UK London}
}

@article{bottcher2022ai,
  title={AI Pontryagin or how artificial neural networks learn to control dynamical systems},
  author={B{\"o}ttcher, Lucas and Antulov-Fantulin, Nino and Asikis, Thomas},
  journal={Nature communications},
  volume={13},
  number={1},
  pages={333},
  year={2022},
  publisher={Nature Publishing Group UK London}
}


%% file: references.bib
@Misc{amsmath,
  author =	 {{American Mathematical Society}},
  title =	 {User's Guide for the \texttt{amsmath} Package
                  (Version 2.0)},
  url =		 {ftp://ftp.ams.org/pub/tex/doc/amsmath/amsldoc.pdf},
  urldate =	 {2015-07-30},
  year =	 2002}

@article{zhang2025time,
  title={A Time-Dependent Symplectic Network for Nonconvex Path Planning Problems with Linear and Nonlinear Dynamics},
  author={Zhang, Zhen and Wang, Chenye and Liu, Shanqing and Darbon, J{\'e}r{\^o}me and Karniadakis, George Em},
  journal={SIAM Journal on Scientific Computing},
  volume={47},
  number={4},
  pages={C769--C794},
  year={2025},
  publisher={SIAM}
}

@book{bardi2008optimal,
  author    = {Martino Bardi and Italo Capuzzo-Dolcetta},
  title     = {Optimal Control and Viscosity Solutions of {Hamilton-Jacobi-Bellman} Equations},
  series    = {Modern Birkh{\"a}user Classics},
  publisher = {Birkh{\"a}user Boston},
  address   = {Boston, MA},
  year      = {2008},
  isbn      = {9780817647544}
}

@article{Crandall1984TwoAO,
  title={Two approximations of solutions of {Hamilton-Jacobi} equations},
  author={M. Crandall and P. Lions},
  journal={Mathematics of Computation},
  year={1984},
  volume={43},
  pages={1-19}
}

@book{falcone2014semi,
  author    = {Maurizio Falcone and Roberto Ferretti},
  title     = {Semi-{L}agrangian Approximation Schemes for Linear and {Hamilton--Jacobi} Equations},
  publisher = {Society for Industrial and Applied Mathematics},
  address   = {Philadelphia, PA},
  year      = {2014},
  isbn      = {9781611973044},
  doi       = {10.1137/1.9781611973051}
}

@article{fleming2000max,
  title={A Max-Plus-Based Algorithm for a {Hamilton--Jacobi--Bellman} Equation of Nonlinear Filtering},
  author={Fleming, Wendell H and McEneaney, William M},
  journal={SIAM Journal on Control and Optimization},
  volume={38},
  number={3},
  pages={683--710},
  year={2000},
  publisher={SIAM}
}

@article{akian2008max,
  title={The max-plus finite element method for solving deterministic optimal control problems: basic properties and convergence analysis},
  author={Akian, Marianne and Gaubert, St{\'e}phane and Lakhoua, Asma},
  journal={SIAM Journal on Control and Optimization},
  volume={47},
  number={2},
  pages={817--848},
  year={2008},
  publisher={SIAM}
}

@article{Mc2007,
  title = {A {{Curse}}-of-{{Dimensionality}}-{{Free Numerical Method}} for {{Solution}} of {{Certain HJB PDEs}}},
  author = {McEneaney, William M.},
  year = {2007},
  month = jan,
  volume = {46},
  pages = {1239--1276},
  issn = {0363-0129, 1095-7138},
  doi = {10.1137/040610830},
  journal = {SIAM Journal on Control and Optimization},
  language = {en},
  number = {4}
}

@article{kang2021algorithms,
  title={Algorithms of data generation for deep learning and feedback design: A survey},
  author={Kang, Wei and Gong, Qi and Nakamura-Zimmerer, Tenavi and Fahroo, Fariba},
  journal={Physica D: Nonlinear Phenomena},
  volume={425},
  pages={132955},
  year={2021},
  publisher={Elsevier}
}

@article{bokanowski2023neural,
  title={Neural networks for first order HJB equations and application to front propagation with obstacle terms},
  author={Bokanowski, Olivier and Prost, Averil and Warin, Xavier},
  journal={Partial Differential Equations and Applications},
  volume={4},
  number={5},
  pages={45},
  year={2023},
  publisher={Springer}
}

@article{raymond1998pontryagin,
  title={Pontryagin's principle for state-constrained control problems governed by parabolic equations with unbounded controls},
  author={Raymond, Jean-Pierre and Zidani, Hasnaa},
  journal={SIAM Journal on Control and Optimization},
  volume={36},
  number={6},
  pages={1853--1879},
  year={1998},
  publisher={SIAM}
}

@article{raymond1999pontryagin,
  title={Pontryagin's principle for time-optimal problems},
  author={Raymond, Jean-Pierre and Zidani, H},
  journal={Journal of Optimization Theory and Applications},
  volume={101},
  number={2},
  pages={375--402},
  year={1999},
  publisher={Springer}
}

@article{alla2019efficient,
  title={An efficient DP algorithm on a tree-structure for finite horizon optimal control problems},
  author={Alla, Alessandro and Falcone, Maurizio and Saluzzi, Luca},
  journal={SIAM Journal on Scientific Computing},
  volume={41},
  number={4},
  pages={A2384--A2406},
  year={2019},
  publisher={SIAM}
}

@article{bokanowski2022optimistic,
  title={Optimistic planning algorithms for state-constrained optimal control problems},
  author={Bokanowski, Olivier and Gammoudi, Nidhal and Zidani, Hasnaa},
  journal={Computers \& Mathematics with Applications},
  volume={109},
  pages={158--179},
  year={2022},
  publisher={Elsevier}
}

@article{akian2024multilevel,
  title={A Multilevel Fast Marching Method for the Minimum Time Problem},
  author={Akian, Marianne and Gaubert, St{\'e}phane and Liu, Shanqing},
  journal={SIAM Journal on Control and Optimization},
  volume={62},
  number={6},
  pages={2963--2991},
  year={2024},
  publisher={SIAM}
}

@article{akian2023adaptive,
  title={An adaptive multi-level max-plus method for deterministic optimal control problems},
  author={Akian, Marianne and Gaubert, St{\'e}phane and Liu, Shanqing},
  journal={IFAC-PapersOnLine},
  volume={56},
  number={2},
  pages={7448--7455},
  year={2023},
  publisher={Elsevier}
}

@article{raissi2019physics,
  title={Physics-informed neural networks: A deep learning framework for solving forward and inverse problems involving nonlinear partial differential equations},
  author={Raissi, Maziar and Perdikaris, Paris and Karniadakis, George E},
  journal={Journal of Computational physics},
  volume={378},
  pages={686--707},
  year={2019},
  publisher={Elsevier}
}

@article{karniadakis2021physics,
  title={Physics-informed machine learning},
  author={Karniadakis, George Em and Kevrekidis, Ioannis G and Lu, Lu and Perdikaris, Paris and Wang, Sifan and Yang, Liu},
  journal={Nature Reviews Physics},
  volume={3},
  number={6},
  pages={422--440},
  year={2021},
  publisher={Nature Publishing Group UK London}
}

@article{cuomo2022scientific,
  title={Scientific machine learning through physics--informed neural networks: Where we are and what’s next},
  author={Cuomo, Salvatore and Di Cola, Vincenzo Schiano and Giampaolo, Fabio and Rozza, Gianluigi and Raissi, Maziar and Piccialli, Francesco},
  journal={Journal of Scientific Computing},
  volume={92},
  number={3},
  pages={88},
  year={2022},
  publisher={Springer}
}

@article{toscano2025pinns,
  title={From pinns to pikans: Recent advances in physics-informed machine learning},
  author={Toscano, Juan Diego and Oommen, Vivek and Varghese, Alan John and Zou, Zongren and Ahmadi Daryakenari, Nazanin and Wu, Chenxi and Karniadakis, George Em},
  journal={Machine Learning for Computational Science and Engineering},
  volume={1},
  number={1},
  pages={1--43},
  year={2025},
  publisher={Springer}
}

@article{lee2024automatic,
  title={Automatic discovery of optimal meta-solvers via multi-objective optimization},
  author={Lee, Youngkyu and Liu, Shanqing and Darbon, Jerome and Karniadakis, George Em},
  journal={arXiv preprint arXiv:2412.00063},
  year={2024}
}

@article{lee2025automatic,
  title={Automatic discovery of optimal meta-solvers for time-dependent nonlinear PDEs},
  author={Lee, Youngkyu and Liu, Shanqing and Darbon, Jerome and Karniadakis, George Em},
  journal={arXiv preprint arXiv:2507.00278},
  year={2025}
}
